\documentclass{article}
\usepackage{graphicx} 
\usepackage{subcaption}
\usepackage{graphicx}
\usepackage[square,numbers]{natbib}
\usepackage[colorlinks=true, linkcolor=blue, citecolor=blue, urlcolor=blue]{hyperref}

\usepackage{amsmath}
\usepackage{amssymb}
\usepackage{xcolor}
\usepackage{bbm}

\usepackage{blindtext}
\usepackage{subfiles}

\usepackage{todonotes}

\usepackage{titlesec}
\titleformat*{\section}{\large\bfseries}
\titleformat*{\subsection}{\normalsize\bfseries}

\usepackage{geometry}
\title{Optimal Reactive Operation of General Topology Supply Chain and Manufacturing Networks under Disruptions}

\begin{document}

\maketitle

\begin{center}
Daniel Ovalle$^{{a}}$, Joshua L. Pulsipher$^{{b}}$, Yixin Ye$^{{c,}}$\footnotemark, Kyle Harshbarger$^{{d,}}$\footnotemark, Scott Bury$^{{e}}$, Carl D. Laird$^{{a}}$,  \\
and Ignacio E. Grossmann$^{{a,}}$\footnotemark\\

\vspace{0.10in}

$^{{a}}$ Department of Chemical Engineering, Carnegie Mellon University, Pittsburgh, PA
\vspace{0.10in}

$^{{b}}$ Department of Chemical Engineering, University of Waterloo, Waterloo, ON N2L 3G1, Canada
\vspace{0.10in}

$^{{c}}$ The Dow Chemical Company, Core R\&D, Lake Jackson, TX
\vspace{0.10in}

$^{{d}}$ The Dow Chemical Company, Supply Chain Innovation, Midland, MI

\vspace{0.10in}

$^{{e}}$ The Dow Chemical Company, Core R\&D, Midland, MI
\end{center}
\footnotetext[1]{\parindent=0cm \small{Currently at DoorDash Inc.}}
\footnotetext[2]{\parindent=0cm \small{Currently at Huntsman Corporation.}}
\footnotetext[3]{\parindent=0cm \small{Corresponding author. Email: {\tt grossmann@cmu.edu}}}

\section*{Abstract}

Supply and manufacturing networks in the chemical industry involve diverse processing steps across different locations, rendering their operation vulnerable to disruptions from unplanned events.
Optimal responses should consider factors such as product allocation, delayed shipments, and price renegotiation , among other factors. 
In such context, we propose a multiperiod mixed-integer linear programming model that integrates production, scheduling, shipping, and order management to minimize the financial impact of such disruptions. The model accommodates arbitrary supply chain topologies and incorporates various disruption scenarios, offering adaptability to real-world complexities.
A case study from the chemical industry demonstrates the scalability of the model under finer time discretization and explores the influence of disruption types and order management costs on optimal schedules. This approach provides a tractable, adaptable framework for developing responsive operational plans in supply chain and manufacturing networks under uncertainty.


\section{Introduction} \label{sec:motivation}

\begin{figure}[!htb]
	\includegraphics[width=0.85\textwidth]{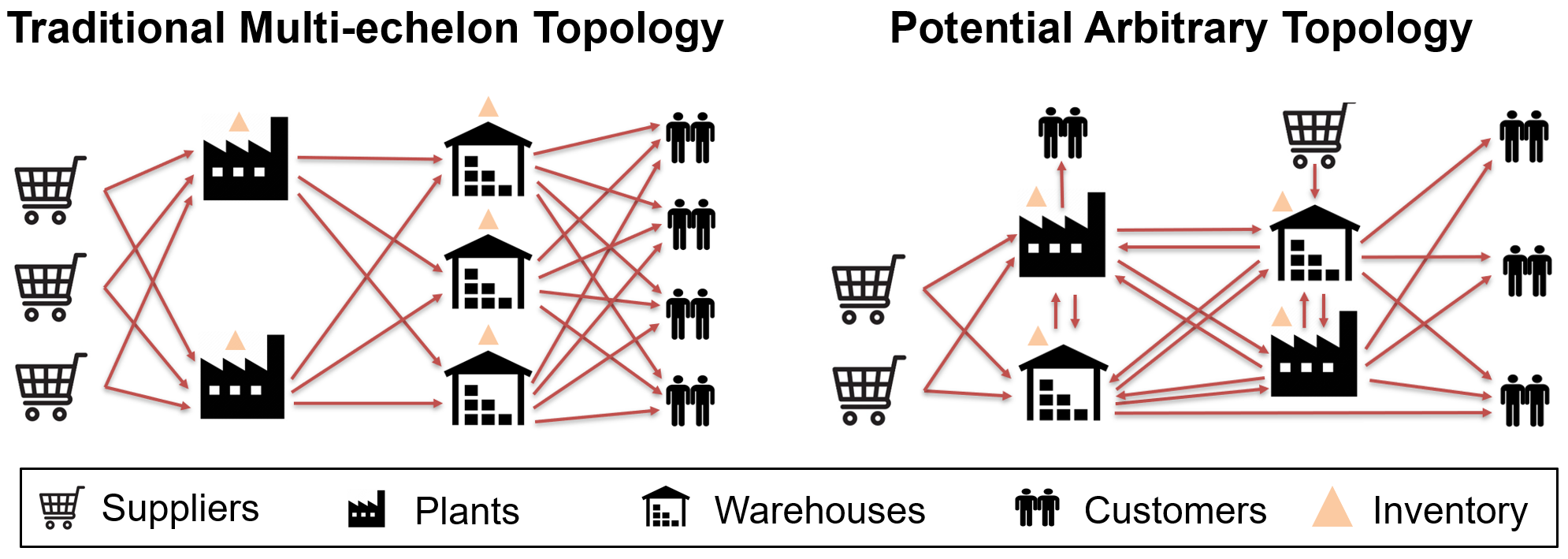}
	\centering
	\caption{Comparison of a traditional feed forward multi-echelon topology an a potential arbitrary topology with similar number of components.}
	\label{fig:topologies}
\end{figure}

Real-world events highlight the importance of effective supply chain disruption management, as exemplified by the COVID-19 pandemic, which introduced unprecedented challenges including facility shutdowns, labor shortages, and transportation constraints \citep{moosavi2022supply,ivanov2020predicting}.
More specifically, supply chain networks of high value-added chemicals are particularly vulnerable to disruption. 
These products often require several raw materials and multiple manufacturing steps in series and/or parallel at various geographical locations. 
Hence, supply chain networks from the chemical industry often deviate from traditional  multi-echelon network topologies where material flows in a feed-forward manner \citep{snyder2019fundamentals}. 
More specifically, these networks have arbitrary network topology, meaning that no simplifying assumptions are made about the structure or connectivity of its components \citep{you2009risk}.
As shown in Figure \ref{fig:topologies}, in a traditional multi-echelon supply chain, material is shipped and transformed sequentially, and the permissible connections among node types are defined. 
In contrast, in arbitrary topologies the notion of an echelon becomes obscured as material flow is no longer unidirectional and can take various paths, including looping through multiple facilities.
While material generally flows ``downstream," it may be reused at multiple stages, particularly in scenarios like sustaining a reaction.
Moreover, plant and warehouse nodes are allowed to have both internal and external demand introducing a significant level of complexity that is not commonly addressed in the supply chain literature \cite{snyder2019fundamentals}.
This is given by the fact that in these complex supply and manufacturing networks the demand cannot be easily aggregated or centralized.
Furthermore, the demand for such products is significantly smaller and more volatile than commodity products. 
Accordingly, the profit margins are higher, which gives rise to an even higher loss rate per disruption when compared to supply chains from other industries.
All of these aspects increase both the challenges and the importance of managing disruptions in high-value chemical supply chains.

From the perspective of the customer, supply chain disruptions and interruptions can lead to delayed deliveries, increases in prices, or impacts on the product quality.
From the perspective of the company, they can lead to significant economic losses, including increased costs to rush shipment to meet customer orders, penalties from delayed deliveries, or even loss of revenue due to failure in fulfilling existing orders.
Furthermore, customer satisfaction levels can be significantly impacted when an order is delivered late or canceled. 
Therefore, effective optimization of decisions along the supply chain is needed to mitigate the economic and operational impacts of disruptions and interruptions. 
Some examples of these decisions include rerouting current inventories, increasing production rate at non-disrupted locations, buying intermediate or finished products from competitors or third parties, and managing existing order delivery dates.
In practice, response decisions are often made locally and empirically leading to sub-optimal strategies. Furthermore, cascading impacts of decentralized decisions can lead to poor system economics, and even situations heuristic policies are not possible to implement.

To address the problem of supply chains under disruptions, a significant amount of literature focuses on supply chain resilience which refers to the ability to adapt to unexpected events and quickly recover the desired level of operations \citep{elleuch2016resilience,kamalahmadi2016review}. 
Researchers have introduced metrics to measure resilience, considering aspects like network information \citep{cardoso2015resilience,adenso2012impact}, the handled products \citep{azevedo2016larg}, and overall performance \citep{munoz2015quantification,dixit2016performance}.
For a more detailed exploration of supply chain resilience, interested readers can refer to the comprehensive reviews by \citet{katsaliaki2021supply} and \citet{ribeiro2018supply}. 

In the area of supply chain resilience, significant attention is given to the facility location problem, which deals with strategic planning for the network design while considering uncertainty associated with demand and disruptions \citep{owen1998strategic, snyder2006facility, snyder2005reliability}. 
This problem balances the trade-off between centralizing the demand on specific nodes to overcome individual site variability \citep{tsiakis2001design, you2008mixed} with the implicit vulnerability that arises at the aggregated critical nodes \citep{snyder2006facility}.
To address this problem, authors have suggested various strategies, including considering site-dependent distribution probabilities \citep{lim2010facility, cui2010reliable}, using approximation algorithms \citep{shen2011reliable, li2010continuum}, implementing tailored cuts and bounds \citep{garcia2014design}, and employing stochastic programming approaches using sample average approximation \citep{santoso2005stochastic, klibi2012scenario}. 
For a more comprehensive review of robust supply chain design under uncertainty, readers can refer to \citet{klibi2010design}. 
While the aforementioned approaches effectively address resilient supply chain network design, gaps remain in the literature regarding short-term operational planning functioning within the restrictions of the physical supply chain design.
There is a need to understand (and optimize) how the operation of the supply chain system can effectively react to given disruptions within a shorter time frame (e.g., daily or hourly) instead of solving a long-term network design problem.

In the operational supply chain literature from the Operations Research (OR) community, several authors derive inventory policies aimed at mitigating the impacts of uncertainty and disruptions, as highlighted by \citet{glasserman1995sensitivity}. 
Utilizing knowledge about network topology proves beneficial in designing and operating supply chains under disruptions. 
When the network exhibits a specific multi-echelon structure, assumptions can be made about its behavior, allowing for development of tailored algorithms or heuristics based on graph and queuing theory \citep{snyder2019fundamentals, snyder2016or}.
Various heuristics and exact algorithms have been developed to determine optimal replenishment and allocation policies under uncertainty for spanning serial networks \citep{chen1994evaluating, clark1960optimal, shang2003newsvendor}, single supplier three-component networks \citep{kim2005adaptive}, assembling networks \citep{rosling1989optimal}, multi-location networks \citep{schmitt2015centralization}, and distribution networks \citep{sherbrooke1968metric, rong2017heuristics}. 
Although these approaches can be optimal and computationally inexpensive, they are not applicable to the case of interest due to two main reasons.  
First, supply chains in the chemical industry frequently have arbitrary topologies, and the assumptions made when developing these methods are not applicable. 
Second, these approaches often yield fixed policies rather than the optimal dynamic response to recover from a specific disruption, leaving a gap in the literature that still needs to be addressed. 
Furthermore, traditional OR rely on distributional assumptions about the demand during lead time which might not always hold in practice \citep{eruguz2016comprehensive, achkar2024efficient}. 


Chemical industry supply chains have been widely addressed by the Process Systems Engineering community \citep{shah2005process, pistikopoulos2021process} as they play a key role in efforts towards enterprise-wide optimization \citep{grossmann2012advances}. 
While the majority of contributions from this community focus on network design, several authors have addressed the operational dynamic response to disruptions as described in the review by \citet{garcia2015supply}.
The operational response was first addressed in \citet{perea2003model}, studying multiperiod supply chain reactions and exploring how closed-loop control techniques can mitigate the bullwhip effect in a multi-echelon supply chain. 
More recently, \citet{badejo2022mathematical} proposed a multiperiod mathematical model coupled with a rolling horizon approach to tackle supply chain disruptions. 
However, there are several key features that this work did not address. 
First, the coarse discretization is not able to capture more detailed dynamics of the supply chain system.
Moreover, this work does not consider important order management decisions, such as late delivery and order cancellation in their response portfolio.

In summary, there is a  need for models capable of generating an optimal dynamic response for operating a supply chain network with an arbitrary topology following a disruption.
This response should include optimal schedules in a actionable time-discretization for production, material acquisition, shipment, and inventory levels. It should also address order management decisions, covering both delivery dates and cancellations of existing orders.
To bridge the aforementioned gaps, in this paper we develop a multiperiod mixed-integer linear programming (MILP) model that determines an optimal disruption-reactive schedule, throughout a specified operation horizon, for a multi-material supply chain and manufacturing network of arbitrary topology. 
This model considers actionable time discretizations that can capture shipments with different transportation durations.
Furthermore, the model is able to account for late deliveries that may accumulate after a disruption, and it is able to provide a schedule on how to optimally deliver those backorders. 
The model also includes binary variables for order cancellation, used when demands cannot be met post-disruption and late delivery is not cost-effective. This strategic cancellation approach allows scarce resources to be reallocated to fulfill critical or high-value orders and enables supply chain managers to notify customers in advance, facilitating opportunities for renegotiation or rescheduling.
To achieve this objective, we use Generalized Disjunctive Programming (GDP) to express discrete decisions at a higher level. 
The resulting disjunctions are then reformulated directly using MILP constraints. For further background on GDP, readers are referred to \citet{raman1994modelling}, and \citet{grossmann2013systematic}.

The remainder of the paper is as follows: Section  \ref{sec:statement} defines the problem, outlining assumptions, decisions, and primary objectives. 
In Section \ref{sec:base_model}, we first introduce the base MILP model and later relax some assumptions to capture more complex supply chain dynamics. 
Section \ref{sec:disruption_modeling} defines various disruption types and demonstrates their natural integration into the model. 
Section \ref{sec:case_study} presents a case study inspired by The Dow Chemical Company's silicone rubber business. 
Finally, Section \ref{sec:results} discusses computational results, and Section \ref{sec:conclusions} provides concluding remarks and suggests potential future research directions.

\section{Problem Statement} \label{sec:statement}


In this work we consider a multi-material supply chain and manufacturing network where a disruption has occurred. 
The network has an arbitrary topology consisting of supplier, plant, warehouse, and customer nodes (see Figure \ref{fig:topologies}). 
We allow for multiple transportation modes, which are represented as different connections between the same nodes.
An operational time horizon is defined during which scheduled customer orders are known, and the reactive operation occurs. 
The shipment times and costs are known for every arc. 
Production coefficients, indicating the relationship between plant-produced recipes and material quantities, are given. 
Capacities for production, shipment, and inventory are specified. 
The material selling price and costs related to material acquisition, production, and inventory holding are also known.

There are several operational variables that need to be determined in the response to disruption. 
Optimal schedules for shipment amounts, routes, plant production, and acquisition need to be constructed for each material at every node. 
By determining these variables alongside existing stock, inventory levels are specified throughout the operation.
Furthermore, given that once a disruption happens it might not always be feasible to deliver orders according to the original schedule, order management decisions need to be made. 
This means determining which orders are fulfilled on time, which will be delayed (and by how much), and which are to be canceled.

The primary goal of this problem is to derive an optimal operational response to the disruption, that aims to maximize the overall profit of the operation throughout the entire horizon while still accounting for customer satisfaction. 
At the end of the time horizon, the obtained response should induce recovery by bringing the system as close as possible to its pre-disruption state. 
This operational strategy involves balancing the existing trade-off between maintaining high customer satisfaction, characterized by delivering orders on time and avoiding cancellations, and obtaining a cost-effective operation during and after the disruption. 

Certain assumptions are made in the problem formulation. 
The production is assumed to be continuous, thus rigorous batch scheduling models such as state-task \citep{maravelias2003new} and resource-task network models \citep{kondili1993general} are not considered. 
Moreover, orders from the same customer and involving the same material within a given time period, are aggregated and treated as a singular entity. 
Additionally, transportation is considered continuous, nevertheless, any discrete unavailability of transportation resources can be effectively modeled as a disruption within the framework.

\section{Proposed Base Model} \label{sec:base_model}
In this section, we present the multiperiod mixed-integer linear programming (MILP) formulation for the optimal integrated response given a disruption. 
This model considers a set of different materials $\mathcal{M}$, as well as suppliers $\mathcal{S}$, plants $\mathcal{P}$, warehouses $\mathcal{W}$, and customers $\mathcal{C}$ that are connected through a set of arcs $\mathcal{A}$ that operate throughout a discrete time horizon $\mathcal{T}$.

The base model operates under several further assumptions.
The model assumes transportation costs are linear with respect to the amount moved and have no minimum shipping amount. Furthermore, fixed costs associated with transportation are ignored.
The base model does not consider service level agreements with suppliers, meaning there is no mandated minimum quantity for procurement from suppliers. 
To facilitate recovery post-disruption, the base model aligns end-of-operation inventories with their pre-disruption levels, i.e., the beginning of the operation. 
Finally, the model imposes a strict lower bound on the buffer stock value, preventing inventories from falling below this threshold. 
All of these assumptions are relaxed in Section \ref{sec:extensions} where different extensions to this base model are discussed. A comprehensive summary of the nomenclature for all the sets, parameters and variables can be found in the Appendix.

\subsection{Material-Node and Material-Arc Set Encoding}

\begin{figure}[!htb]
	\includegraphics[width=0.9\textwidth]{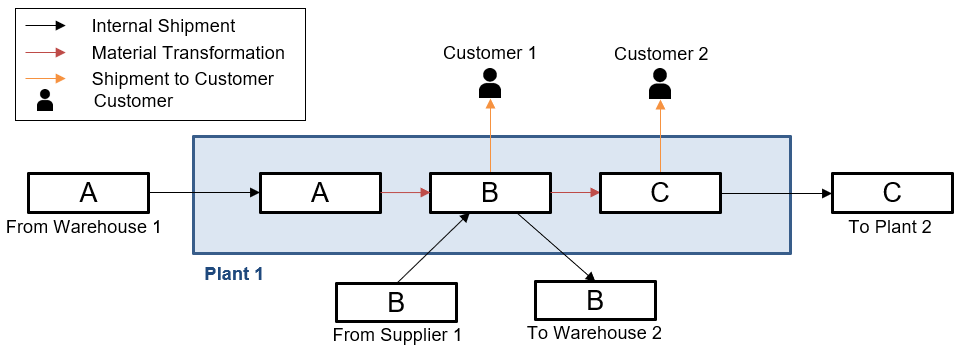}
	\centering
	\caption{Illustrative example of material transformation at Plant 1.}
	\label{fig:material_flexibility}
\end{figure}

Flexibility in managing various materials is crucial in supply chain optimization with arbitrary network topologies. 
Conventional distinctions between raw materials, intermediate products, and final products can become ambiguous, as the same material may serve different functions based on its location within the supply chain.
For example, consider a supply chain network with a manufacturing facility, Plant 1, as depicted in Figure \ref{fig:material_flexibility}.

In this plant, material B can be viewed either as a raw material purchased from Supplier 1, or as an intermediate product from the internal reaction. 
Furthermore, it can be considered as a product given
external demand from Customer 1. 
Figure \ref{fig:material_flexibility} emphasizes the need of adopting a flexible material handling approach, avoiding rigid classifications. 

In response to this, the proposed model employs a single aggregated set $\mathcal{M}$ containing all materials.
This allows for material tracking irrespective of its specific location or  role within the network. Moreover, to avoid creation of unnecessary variables (e.g., materials that are not sold by a supplier), we define subsets of materials existing at a given node $n$ as $\mathcal{M}_n \ \forall \ n\in \mathcal{N}$.
Similarly, not every material has to traverse every arc, leading to the definition of the analogous subset $\mathcal{M}_a \ \forall \ a\in \mathcal{A}$, as all the species capable of traveling through arc $a$. 
Also, we incorporate the transportation method directly into the arc definition, where each arc is characterized by origin, destination, material, and transportation method. 
This approach, together with the set of materials $\mathcal{M}_a$, provides the modeler with flexibility, allowing the designation of materials that can exclusively traverse a route when shipped via a particular transportation mode.

\subsection{Balance Constraints}

The parameter $\tau_{mat}$ represents the shipping time (duration) for material $m$ entering arc $a$ at time period $t$.
Also, $F^{In}_{mat}$ represents the amount of material $m$ entering arc $a$ during a given time period $t$, while $F^{Out}_{mat}$ represents the same flow emerging from the arc after $\tau_{mat}$ time periods as:
\begin{equation}
    F^{In}_{mat} = F^{Out}_{ma\{t+\tau_{mat}\}} \quad \forall \ a \in \mathcal{A}, \ m \in \mathcal{M}_a, \ t \in \mathcal{T} : t+\tau_{mat} \leq |\mathcal{T}|. \label{eq:delay}
\end{equation}
This also allows us to model transportation delays as disruptions within the framework.
Note also that extra constraints are included to disable flow variables that will lead to times outside the time horizon: 
\begin{subequations}
\label{eq:delays}
\begin{equation}
    F^{In}_{mat} = 0 \quad \forall \ a \in \mathcal{A}, \ m \in \mathcal{M}_a, \ t \in \mathcal{T} : t+\tau_{mat} > |\mathcal{T}| \label{eq:delay1}
\end{equation}
\begin{equation}
    F^{Out}_{mat} = 0 \quad \forall \ a \in \mathcal{A}, \ m \in \mathcal{M}_a, \ t \in \mathcal{T} : t-\tau_{mat} < 0 \label{eq:delay2}
\end{equation}
\end{subequations}

The warehouse inventory ($I$) balance is shown below where $F^{Out}$ variables include flows going out of an upstream arc into the inventory, and $F^{In}$ refer to flows out of the inventory into a downstream arc.

\begin{equation}
\label{eq:warehouse_bal}
    I_{mwt} = I_{mw\{t-1\}} + \sum_{a \in \mathcal{A}^{In}_w} F^{Out}_{mat} - \sum_{a \in \mathcal{A}^{Out}_w} F^{In}_{mat} \quad \forall \ w \in \mathcal{W}, \ m \in \mathcal{M}_w, \ t \in \mathcal{T} \backslash \{0\}.
\end{equation}

Plants are allowed to store inventory and utilize existing stock in their reactive response.
This flexibility allows for policies, such as producing with available inventory to alleviate upstream shortages, or accumulating production in the early stages of a reaction for later depletion. 
Consequently, an inventory balance similar to Equation \eqref{eq:warehouse_bal} is used along with an additional consumption/generation term to account for production as:


\begin{equation}
    I_{mpt} = I_{mp\{t-1\}} + \sum_{a \in \mathcal{A}^{In}_p} F^{Out}_{mat}  - \sum_{a \in \mathcal{A}^{Out}_p} F^{In}_{mat}
  + \sum_{r \in \mathcal{R}_p} \phi_{rm}P_{prt}  \quad \forall \ p \in \mathcal{P}, \ m \in \mathcal{M}_p, \ t \in \mathcal{T} \backslash \{0\}, \label{eq:plant_bal_simple}
\end{equation}
where production balance ($P$) is modeled using a recipe-based approach, and the set $\mathcal{R}_p$ denotes the possible recipes at plant $p$.
In this method, $P$ represents the amount of a recipe produced. To map from recipe to material, the mass-based bill of materials $\phi$ is employed, indicating whether a material is consumed (-) or produced (+).
 This method allows for the modeling of different recipes by tracking the net production or consumption of each material.
In the base model production is assumed to occur instantaneously per time period, however, an extension to time-dependent production is discussed in Section \ref{sec:production_time}.
 

Inputs and outputs are calculated in the network in terms the flows leaving supplier nodes and the flows entering demand nodes, respectively. 
This allows every plant and warehouse to either receive external supply, or fulfill external demand, as long as the corresponding arc to the supplier or customer exists. 
The satisfied demand ($D$) is then computed considering all arcs entering a customer $c$ ($A^{In}_c$) through the equation:

 \begin{equation}
    D_{mct} = \sum_{a \in \mathcal{A}^{In}_c} F^{Out}_{mat} \quad \forall \ c \in \mathcal{C}, \ m \in \mathcal{M}_c, \ t \in \mathcal{T}. \label{eq:demand}
\end{equation}
Here, $D_{mct}$ determines the customer demand amount that is satisfied at time $t$, and may or may not be the same as the desired amount ordered. This is further discussed in Section \ref{sec:order_management}.
Analogously, the amount to buy ($B$) is computed by accounting for all arcs leaving supplier $s$ ($A^{Out}_s$) as inidicated below:
\begin{equation}
      B_{mst} = \sum_{a \in \mathcal{A}^{Out}_s} F^{In}_{mat} \quad \forall \ s \in \mathcal{S}, \ m \in \mathcal{M}_s, \ t \in \mathcal{T}. \label{eq:buy}
\end{equation}

\subsection{Determining Initial and Final States}

The proposed formulation is a reactive model that plans supply chain operations once disruption information becomes available. 
Hence, once there is awareness of a disruption, the initial state of the system, characterized by the inventory levels in both plants and warehouses at the time the unplanned event takes place ($I^0$), is fixed as:
\begin{align}
     I_{mn0}= I_{mn}^0  \quad \forall \ n \in \mathcal{P} \cup \mathcal{W}, \ m \in \mathcal{M}_n.  \label{eq:initial_state}
\end{align}

  Similarly, the initial state of all remaining decision variables is determined by their pre-planned values in the absence of disruption as:

 \begin{subequations}
 \label{eq:fixing_to_zero}
 \begin{equation}
    F^{In}_{ma0} = F^{In,0}_{ma} \quad \forall \ a \in \mathcal{A}, \ m \in \mathcal{M}_a \label{eq:Fin_fix}
\end{equation}
\begin{equation}
      F^{Out}_{ma0} = F^{Out,0}_{ma} \quad \forall \ a \in \mathcal{A}, \ m \in \mathcal{M}_a \label{eq:Fout_fix}
\end{equation}
\begin{equation}
     P_{pr0} =  P_{pr}^0 \quad \forall \ p \in \mathcal{P}, \ r \in \mathcal{R}_p \label{eq:P_fix}
\end{equation}
\begin{equation}
      D_{mc0} = D_{mc}^0  \quad \forall \ c \in \mathcal{C}, \ m \in \mathcal{M}_c \label{eq:D_fix}
\end{equation}
\begin{equation}
      B_{ms0} = B_{ms}^0 \quad \forall \ s \in \mathcal{S}, \ m \in \mathcal{M}_s \label{eq:B_fix}
\end{equation}
\end{subequations}

In practice, this approach can be deployed via a rolling horizon, or a shrinking horizon technique \citep{you2009risk, balasubramanian2004approximation}. 
These strategies not only allow decision-makers to update model parameters with current information but also facilitate the incorporation of new disturbances, providing the robust advantages of feedback \citep{rawlings2017model}.

 Ensuring a notion of recovery post-disruption is important to facilitate a smoother transition to nominal inventory policies after the disruption has been addressed. Defining a sense of normality in a supply chain system prone to constant perturbations can be subject to multiple interpretations. Thus, in this work our objective is to bring the system back to the state it was before the disruption occurred.
In the base model we do this by imposing hard terminal constraints as:
 \begin{align}
    I_{mn|\mathcal{T}|} = I_{mn}^0  \quad \forall \ n \in \mathcal{P} \cup\mathcal{W}, \ m \in \mathcal{M}_n,  \label{eq:final_time}
\end{align}
 which have been used successfully to recover disrupted systems in the literature \citep{chong2009model}. Nevertheless, since this approach may appear overly restrictive to certain modelers, we propose an alternative method for handling the final state in Section \ref{sec:final_inventory}.

If the final state is not handled properly, the optimizer may  drive inventories to zero once the reactive operation concludes due to the existence of holding costs in the objective function \citep{lima2011long}. 
Furthermore, the practical deployment of this model can be perceived as an instantiation of model predictive control, wherein solution profiles (i.e., schedules) are derived at each time step, and these solutions are implemented using a rolling horizon strategy \citep{perea2003model, subramanian2013integration}. 
In such a setup, ensuring stability during the rolling horizon requires proper handling of the end time point, and the use of a hard terminal constraint in Equation \eqref{eq:final_time} constitutes a viable option \citep{rawlings2017model}. 
Furthermore, the following subsection, we introduce additional order management options that ensure feasibility, even with the imposition of a hard terminal constraint.



\subsection{Order Management Dynamics} \label{sec:order_management}
An important innovation of this work is accounting for late deliveries and order cancellations as potential reactive actions, which is key given that after a disruption occurs, adhering to the original order schedule may become infeasible or cost prohibitive.

In this model, orders may be satisfied later than their initially planned delivery date by paying a late delivery penalty. 
To achieve this, unmet demand is accumulated, and it can be progressive fulfilled over time as resources become available in the system. 
With this model, fractions of the original order can be delivered to the customer incrementally until the initially accepted order is entirely fulfilled.
To  incentivize the prompt satisfaction of delayed orders, unmet demand is penalized for each period that it remains undelivered.

\begin{figure}[!htb]
\centering
\begin{subfigure}{0.49\textwidth}
    \includegraphics[width=1.1\textwidth]{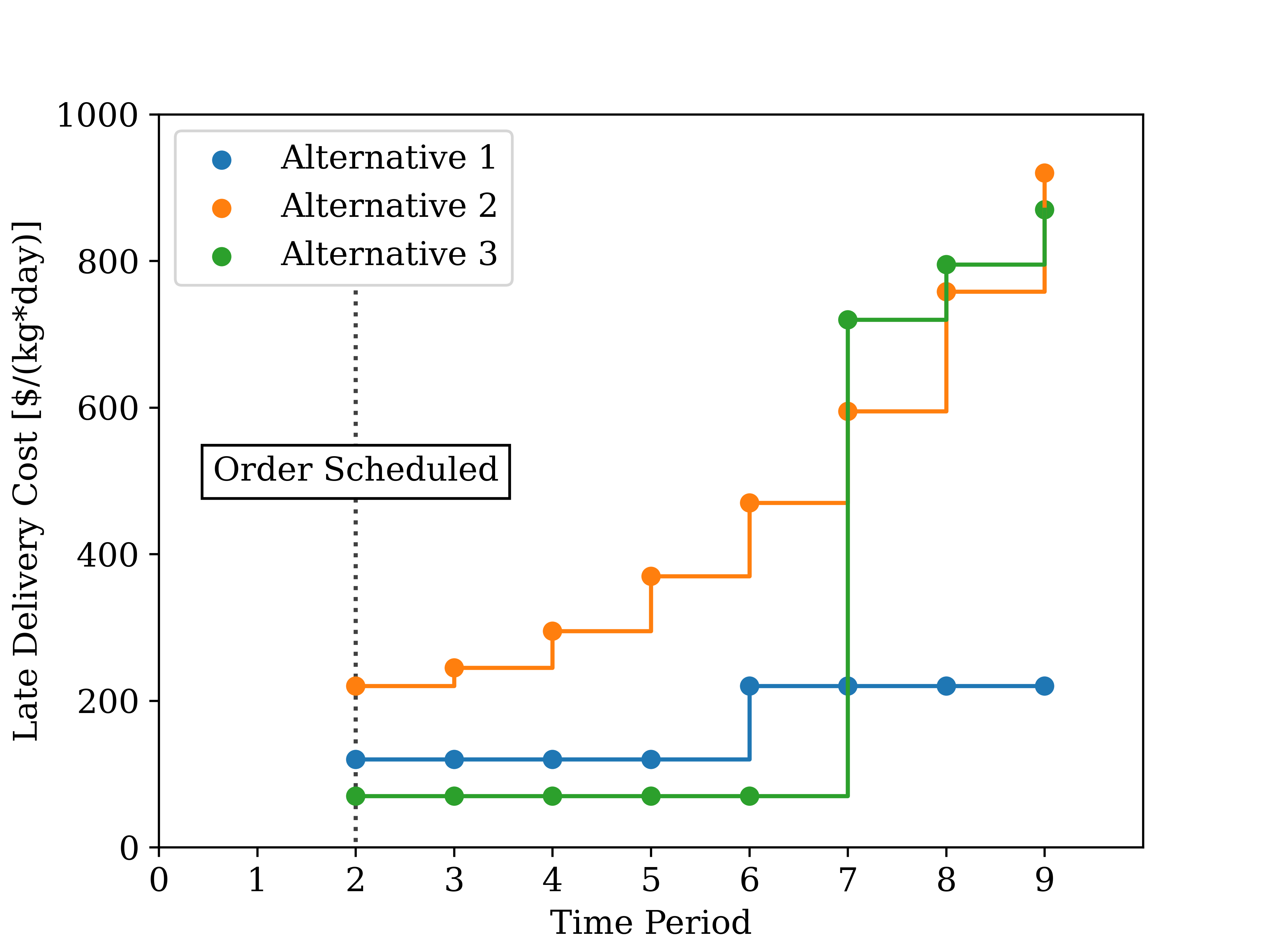}
    \caption{Late delivery cost ($\lambda^U$) profiles.}
    \label{fig:late_delivery_profile}
\end{subfigure}
\hfill
\begin{subfigure}{0.49\textwidth}
    \includegraphics[width=1.1\textwidth]{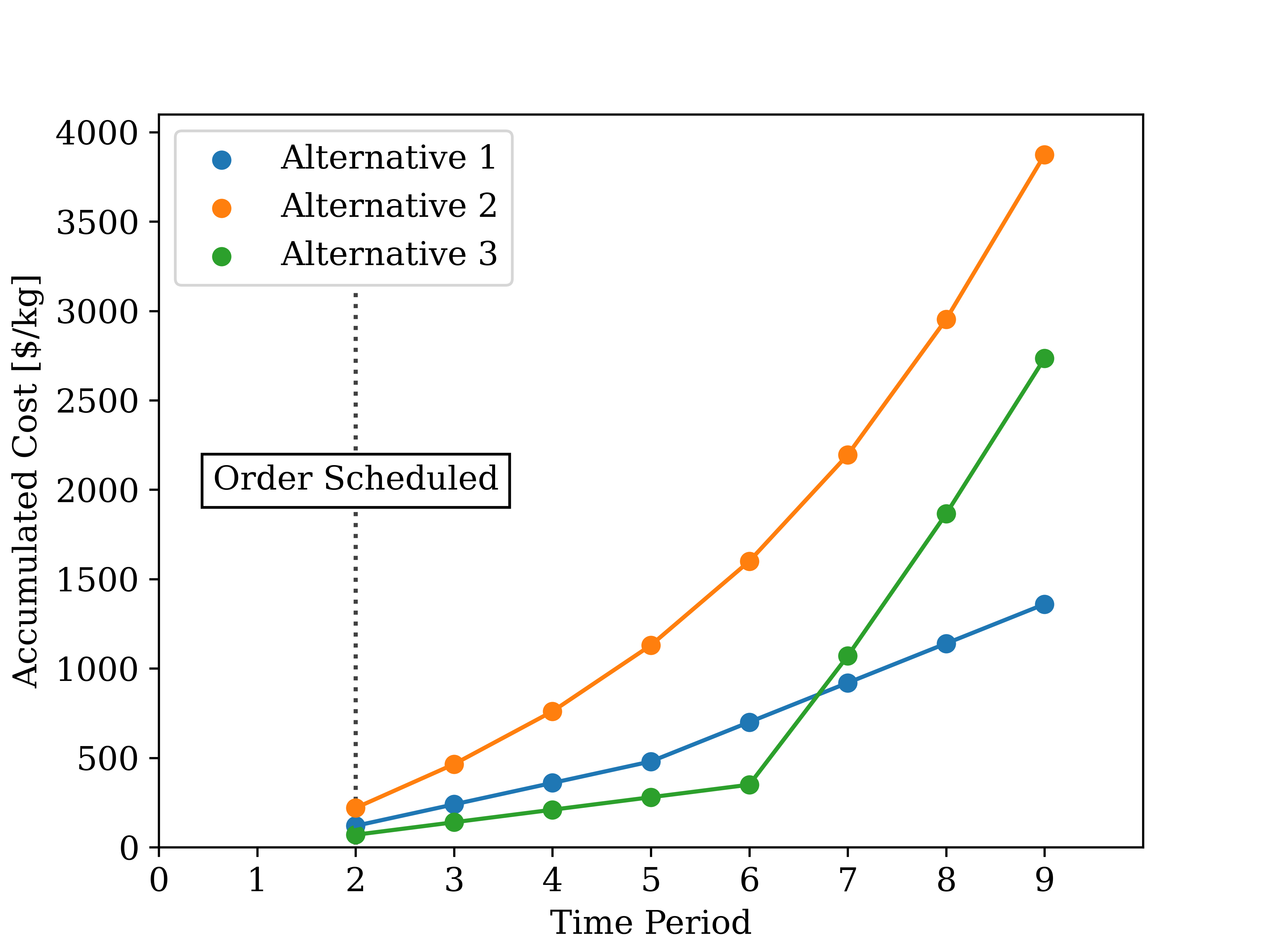}
    \caption{Respective late delivery cost accumulation.}
    \label{fig:late_delivery_accumulation}
\end{subfigure}
        
\caption{Potential late delivery cost profiles and respective cost accumulations to account for different customer priorities.}
\label{fig:late_delivery}
\end{figure}

The reactive response can be customized for individual customers and products by tailoring the late delivery cost profile for each customer and material, as illustrated in Figure \ref{fig:late_delivery_profile}. 
Correspondingly, each profile results in a distinct cost accumulation for units that remain unmet, as depicted in Figure \ref{fig:late_delivery_accumulation}. 
For this, we assume that orders for the same material and customer are sufficiently spaced in time, which is reasonable given in the chemical product supply chain orders for the same product are typically large in size and widely dispersed in time.
This assumption is made to minimize order interference in the profiles. Nevertheless, if this condition does not hold, a constant penalty per customer and material will still provide a sense of prioritization as needed.


The example given in Figure \ref{fig:late_delivery} examines an order scheduled for period 2 for a specific customer and material to illustrate the flexibility in prioritization that can be achieved. 
For simplicity, we show the first order in the schedule post-disruption, hence, no cost profile is required for previous time periods.
Alternative 1 features a constant cost that increases in the second half of the time horizon representing an average customer  where no particular priority is required, maintaining a constant cost with an increased penalty for older unmet demand.
In contrast, Alternative 2 outlines a profile for a higher-priority customer featuring a larger initial cost and an exponentially increasing profile, incentivizing a faster delivery than the previous alternative. 
Alternative 3 reflects a scenario where a customer tolerates late deliveries for five days, keeping the cost low during this grace period, but significantly escalating afterward, showcasing how more complex customer dynamics can be integrated into the model by adjusting the cost profile.
Analogously, Figure \ref{fig:late_delivery_accumulation} illustrates how the late delivery cost accumulates over time for a single unit of unmet demand, providing the total penalty cost associated with delivering at each time period for the previous alternatives.

In scenarios where delivering an order significantly past its original due date is unavoidable, it might be worth considering the option of canceling the order rather than delivering late.
Canceling an order has the potential to alleviate the network, enabling the allocation of resources to other orders more efficiently and leading to a more cost-effective response.
Consequently, our model incorporates a feature for order cancellation, allowing the optimizer to make informed decisions about whether to accept an order, and if accepted, how late it will be progressively (or entirely) fulfilled. 
This flexibility enhances the response of the model, allowing it to optimally manage the fulfillment of orders and their potential cancellation in the face of disruptions.
If the model decides to cancel an order, it incurs a fixed penalty reflecting the impact of this cancellation on the goodwill of the company with its customers.

By leveraging the GDP framework, we can formulate a disjunction for the cancellation and its impact on the unmet demand ($U$) balance as follows:
\begin{equation}
    \left[\begin{gathered}
        Y_{mct}\\
       U_{mct} = U_{mc\{t-1\}} - D_{mct}
    \end{gathered}\right] \vee
    \left[\begin{gathered}
        \neg Y_{mct}\\
        U_{mct} = U_{mc\{t-1\}} - D_{mct} + \delta_{mct}
    \end{gathered}\right] \quad \forall \ c \in \mathcal{C},  \ m \in \mathcal{M}_c, \ t \in \mathcal{T}.
    \label{eq:unmet_disjunction}
\end{equation}

This formulation implies that when an order is canceled ($Y{=}\text{True}$), that order amount $\delta$ should not be accounted for in the balance allowing the operation to focus on the  unmet demand accumulated from previously accepted orders. 
Conversely, when the order is accepted (or not canceled, $Y{=}\text{False}$), the amount ordered $\delta$ is included in the unmet balance. 
This allows the model to determine if the order will be delivered on time or late, and if so, how it will be progressively fulfilled.

As mentioned in Section \ref{sec:motivation}, GDP is introduced as a modeling framework to pose the constraints more clearly. Since, the problem is not solved in disjunctive form, we introduce an equivalent MILP reformulation where $y$ is the binary variable associated with $Y$ (meaning $y = \mathbbm{1}\{Y \}$, where $\mathbbm{1}\{\cdot\}$ is the indicator function) as:
\begin{align}
    U_{mct} = U_{mc\{t-1\}} - D_{mct} + \delta_{mct}(1-y_{mct}) \quad \forall \ c \in \mathcal{C}, \ m \in \mathcal{M}_c, \ t \in \mathcal{T}  \label{eq:unmet_mip}
\end{align}

The resulting expression for the unmet demand is calculated as an inventory balance to account for accumulated unmet demand across the time periods. 
Since unmet demand is nonnegative, the model can choose to cancel a specific order ($y = 1$), relieving the system but incurring a high fixed price ($\lambda^\delta$).
On the other hand, if a given order is not canceled ($y = 0$), the reactive operation tries to fulfill both the given order and any accumulated unmet demand. 
The above expression applies only when an order is present, meaning that $\delta$ is non-zero. If there is no order (i.e., $\delta=0$), the cancellation variable $y$ can be set to one.
The impact of canceled orders on the objective function is discussed in Section \ref{sec:obj_and_bounds}.


\subsection{Inventory bounds}\label{sec:inventory_bounds}
Upper bounds for the inventory levels in plants and warehouses are typically given by the physical capacity allocated to each material within the facility ($I^{U}_{mnt}$). 
However, when materials are not segregated, inventory may be restricted relative to the total facility volume ($I^{U}_{nt}$) as:
\begin{align}
    \sum_{m \in \mathcal{M}_n} I_{mnt} \leq I^{U}_{nt} \quad \forall \ n \in \mathcal{P} \cup \mathcal{W}, \ t \in \mathcal{T}\label{eq:inv_bound_upper}
\end{align}

Determining the lower bound for inventory levels is a decision influenced by the risk aversion of the modeler. 
The riskier approach permits inventories to reach zero, allowing all available buffer stock to be utilized during disruptions if necessary. 
Conversely, the more conservative approach retains a lower bound on the buffer stock ($SS$), based on the argument that buffer stock is primarily intended for hedging against demand uncertainty rather than disruptions. 
Both perspectives can be captured and combined in the following expression:
\begin{equation}
       {I}^L_{mnt} = \alpha_{mnt} SS_{mnt} \leq I_{mnt} \quad \forall \ n \in \mathcal{P}\cup\mathcal{W}, \ m \in \mathcal{M}_n, \ t \in \mathcal{T} \label{eq:inv_bound_lower_hard}
\end{equation}
where $\alpha \in [0,1]$ represents the fraction of the buffer stock $SS$ that the model preserves regardless of the disruption.
The inventory bounds used in the base model are explicitly shown in the next section (see Equation \eqref{eq:I_bounds}) where a capacity is assumed for each material individually and no lower bound is placed to yield the least conservative schedules ($\alpha=0$). Nevertheless, an alternative where this assumption is relaxed is proposed in Section \ref{sec:negative_dev}.

\subsection{Objective Function}\label{sec:obj_and_bounds}

The objective of the model is to maximize profit from the sales ($\lambda^D$), while accounting for several costs such as buying material from suppliers, manufacturing cost at different facilities ($\lambda^P$), shipping cost using different transportation types (road, marine or air-freight) ($\lambda^F$), and holding inventory cost at both plants and warehouses ($\lambda^{I}$). 
Furthermore, this model accounts for cost of delivering late ($\lambda^U$), and for a high fixed cost in the case where an order needs to be canceled ($\lambda^\delta$) as shown below:
\begin{align}
    \max & \sum_{t \in \mathcal{T}} [ \sum_{c \in \mathcal{C}} \sum_{m \in \mathcal{M}_c} (\lambda^D_{mct}D_{mct} - \lambda^U_{mct}U_{mct} - \lambda^\delta_{mct}y_{mct})  - \sum_{a \in \mathcal{A}} \sum_{m \in \mathcal{M}_a} \lambda^F_{mat}F^{Out}_{mat} - \sum_{s \in \mathcal{S}} \sum_{m \in \mathcal{M}_s} \lambda^B_{mst}B_{mst} \nonumber\\ 
  & - \sum_{p \in \mathcal{P}} \sum_{r \in \mathcal{R}_p} \lambda^P_{prt}P_{prt} 
  - \sum_{n \in \mathcal{P} \cup \mathcal{W}} \sum_{m \in \mathcal{M}_n} \lambda^{I}_{mnt}I_{mnt}]  \label{eq:obj}
\end{align}

The values of $\lambda^U$ and $\lambda^\delta$ capture the impact that delays in delivery or order cancellations can have on the reputation of the company and its customer relations.
Determining the exact value for these parameters is a modeling decision that must obtain a balance between financial considerations and customer relationships. 
These penalties have bounds that can be inferred from the system. The upper bound is the net present value of all forecasted sales from the affected customer over 15 to 20 years, representing the worst-case scenario of customer loss due to a canceled order or late delivery. The lower bound is the cost of sourcing the product from a competitor and selling it at cost, covering logistical expenses without affecting goodwill.

Finally, all continuous variables are nonnegative and bounded as shown below: 

 \begin{subequations}
 \label{eq:bounds}
 \begin{equation}
    0 \leq F^{In}_{mat} \leq {F}^U_{mat} \quad \forall \ a \in \mathcal{A}, \ m \in \mathcal{M}_a, \ t \in \mathcal{T} \label{eq:Fin_bounds}
\end{equation}
\begin{equation}
      0 \leq F^{Out}_{mat} \leq {F}^U_{mat} \quad \forall \ a \in \mathcal{A}, \ m \in \mathcal{M}_a, \ t \in \mathcal{T} \label{eq:Fout_bounds}
\end{equation}
\begin{equation}
      0 \leq I_{mnt} \leq {I}^U_{mnt} \quad \forall \ n \in \mathcal{P}\cup\mathcal{W}, \ m \in \mathcal{M}_n, \ t \in \mathcal{T} \label{eq:I_bounds}
\end{equation}
\begin{equation}
      0 \leq P_{prt} \leq {P}^U_{prt} \quad \forall \ p \in \mathcal{P}, \ r \in \mathcal{R}_p, \ t \in \mathcal{T} \label{eq:P_bounds}
\end{equation}
\begin{equation}
      0 \leq D_{mct} \leq {D}^U_{mct} \quad \forall \ c \in \mathcal{C}, \ m \in \mathcal{M}_c, \ t \in \mathcal{T} \label{eq:D_bounds}
\end{equation}
\begin{equation}
      0 \leq B_{mst} \leq {B}^U_{mst} \quad \forall \ s \in \mathcal{S}, \ m \in \mathcal{M}_s, \ t \in \mathcal{T} \label{eq:B_bounds}
\end{equation}
\end{subequations}
where the bounds are time-indexed. Additional considerations on bounding variables to model disruptions will be presented in Section \ref{sec:disruption_modeling}.

\subsection{Modeling Extensions} \label{sec:extensions}

The base model is composed of Equations \eqref{eq:delay}-\eqref{eq:final_time}, \eqref{eq:unmet_mip}-\eqref{eq:bounds} and it provides a reactive supply chain operation solution in response to the disruptions under mild assumptions.
The following subsections illustrate how the base model can be extended to accommodate more flexible responses and dynamics, when some of the initial assumptions are relaxed. Later in this paper (see Section \ref{sec:scalability}), we will discuss the impact that each of the extensions, individually and in combination, has on the computational performance.

\subsubsection{Production Across Time Periods (PATP)}\label{sec:production_time}
The base model assumes time-independent production, meaning that the recipe starts and concludes within the same period as shown in the Equation \eqref{eq:plant_bal_simple}. 
However, there are instances where this assumption may not hold true, either due to processes inside the plant taking longer than the time discretization (e.g., brewing or fermenting), or because the modeler opts for a finer time scale to capture more detailed schedules. 
To model production across time periods (PATP), we follow a similar approach to that used for modeling shipping time delays, as demonstrated in Equations \eqref{eq:delay}-\eqref{eq:delays}. 
We define $P^{In}$ and $P^{Out}$ as the quantity of produced recipe that initiates and concludes (respectively) within the specified period. Now, the time-dependent production can be expressed as:
\begin{subequations}
\label{eq:production_delays}
\begin{equation}
    P^{In}_{prt} = P^{Out}_{pr\{t+\tau^P_{prt}\}} \quad \forall \ p \in \mathcal{P}, \ r \in \mathcal{R}_p, \ t \in \mathcal{T} \label{eq:production_delay}
\end{equation}
\begin{equation}
    P^{In}_{prt} = 0 \quad \forall \ p \in \mathcal{P}, \ r \in \mathcal{R}_p, \ t \in \mathcal{T} : t+\tau^P_{prt} > |\mathcal{T}| \label{eq:production_delay1}
\end{equation}
\begin{equation}
    P^{Out}_{prt} = 0 \quad \forall \ p \in \mathcal{P}, \ r \in \mathcal{R}_p, \ t \in \mathcal{T} : t-\tau^P_{prt} < 0 \label{eq:production_delay2}
\end{equation}
\end{subequations}
where $\tau^P_{prt}$ is the time it takes to execute the recipe when production starts at period $t$.
  Consequently, the plant inventory balance originally expressed in Equation \eqref{eq:plant_bal_simple} can be written in terms of the time-dependent production as:
\begin{equation}
\label{eq:plant_bal_delays}
    I_{mpt} = I_{mp\{t-1\}} + \sum_{a \in \mathcal{A}^{In}_p} \hspace{-1mm} F^{Out}_{mat}  - \sum_{a \in \mathcal{A}^{Out}_p} \hspace{-1mm} F^{In}_{mat} 
   + \hspace{-4mm}  \sum_{r \in \mathcal{R}_p | \phi_{rm}>0} \hspace{-4mm} \phi_{rm}P^{Out}_{prt} + \hspace{-4mm} \sum_{r \in \mathcal{R}_p | \phi_{rm}<0} \hspace{-4mm}\phi_{rm}P^{In}_{prt} \ \ \forall \ p \in \mathcal{P}, \ m \in \mathcal{M}_p, \ t \in \mathcal{T} \backslash \{0\}
\end{equation}

Therefore, the model that supports manufacturing that occurs across time periods is given by Equations \eqref{eq:delay}-\eqref{eq:warehouse_bal}, 
\eqref{eq:demand}-\eqref{eq:final_time}, and 
\eqref{eq:unmet_mip}-\eqref{eq:plant_bal_delays}, where the production variable is redefined.


\subsubsection{Fixed Transportation Cost (FTC) and Minimum Shipping Amount}\label{sec:fixed_transport} 
The base model assumes that  shipping costs are proportional to the quantity of material shipped. However, in many situations, there may be fixed transportation costs (FTC) associated with the use of a particular route. Furthermore, there may also be minimum shipment amount requirements.

To extend the model to account fixed transportation costs and minimum quantities, we define $X$ as a Boolean variable that indicates flow existence according to the following disjunction: 



\begin{equation}
    \left[\begin{gathered}
        X_{mat}\\
       {F}^L_{mat} \leq F^{In}_{mat} \leq {F}^U_{mat}
    \end{gathered}\right] \vee
    \left[\begin{gathered}
        \neg X_{mat}\\
        F^{In}_{mat} = 0
    \end{gathered}\right] \quad \forall \ a \in \mathcal{A}, \ m \in \mathcal{M}_a, \ t \in \mathcal{T}
    \label{eq:fixed_disjunction}
\end{equation}
where flow is activated only when $X=\text{True}$, and it is bounded by both the minimum quantity and the capacity of that arc $[{F^L},{F^U}]$.
Similar to the disjunction shown in Equation \eqref{eq:unmet_disjunction}, we propose an MILP reformulation using the associated binary variable $x$ (where $x = \mathbbm{1}\{X \}$) as:
\begin{align}
    {F}^L_{mat}x_{mat} \leq F^{In}_{mat} \leq {F}^U_{mat}x_{mat} \quad \forall \ a \in \mathcal{A}, \ m \in \mathcal{M}_a, \ t \in \mathcal{T}
    \label{eq:fixed_mip}
\end{align}
allowing the flow to be non-zero only when there is material traversing the arc. 
Also, $x$ can be used to account for fixed transportation costs by subtracting:
\begin{align}
    FixedCost = \sum_{t \in \mathcal{T}}\sum_{a \in \mathcal{A}} \sum_{m \in \mathcal{M}_a} \lambda^{Ffix}_{mat} x_{mat}
    \label{eq:fixed_cost_term}
\end{align}
from the objective function of the base model depicted in Equation \eqref{eq:obj}.

The extension that considers fixed transportation cost and minimum shipping amount is given by Equations \eqref{eq:delay}-\eqref{eq:final_time}, \eqref{eq:unmet_mip}-\eqref{eq:bounds}, and  \eqref{eq:fixed_mip}-\eqref{eq:fixed_cost_term}.
 This extension introduces a significant number of binary variables and big-M constraints (like Equation \eqref{eq:fixed_mip}) which can impact the solution time and the tightness of the linear programming relaxation of the MILP \citep{trespalacios2015improved}. 
The impact of these extensions on the solution time and scalability of the model are shown in Section \ref{sec:scalability}. 

\subsubsection{Service Level Agreements (SLA)}\label{sec:sla}
Service level agreements (SLA) are contracts that may exist with the potential suppliers.
If the supply chain company desires to buy material from the supplier, then a minimum amount may need to be purchased throughout a period of time.

Service level agreements can be modeled in two different ways. 
The first option is to buy more than a minimum permitted amount and can be used when the time discretization (i.e. daily or weekly) of the model is equal or larger than the one in the contract.
In that case, a simple minimum quantity bound is imposed as: 
\begin{equation}
    \left[\begin{gathered}
        W_{mst}\\
       {B}_{mst}^{SLA} \leq {B}_{mst} \leq {B}^U_{mst}
    \end{gathered}\right] \vee
    \left[\begin{gathered}
        \neg W_{mst}\\
       {B}_{mst} = 0
    \end{gathered}\right] \quad \forall \ s \in \mathcal{S}^{SLA}, \ m \in \mathcal{M}_s, \ t \in \mathcal{T}
    \label{eq:sla1_disjunction}
\end{equation}
where $W$ is the Boolean variable that determines if the SLA is executed at time $t$, and the model commits to purchase the required quantities. 
Furthermore, we define $B^{SLA}$ as the minimum purchase amount.

In the second option, where the time partitioning of the model is smaller than the one in the agreement, a time window has to be considered.
Here, the model has to guarantee that the quantity purchased over the SLA time window is within the contract specification. This can be modeled as follows:
\begin{equation}
    \left[\begin{gathered}
        W_{mst}\\
       {B}_{mst}^{SLA} \leq \sum_{t'=t}^{t+\tau^{SLA}_{mst}}{B}_{mst'} \leq {B}^U_{mst}
    \end{gathered}\right] \vee
    \left[\begin{gathered}
        \neg W_{mst}\\
       {B}_{mst} = 0
    \end{gathered}\right] \quad \forall \ s \in \mathcal{S}^{SLA}, \ m \in \mathcal{M}_s, \ t \in \mathcal{T}: t+\tau^{SLA}_{mst} \leq |\mathcal{T}| 
    \label{eq:sla2_disjunction}
\end{equation}
where $\tau^{SLA}$ is the duration of the time window. The equation above assumes that a contract cannot extend beyond the time horizon; thus, the disjunction is only imposed when the SLA can be fulfilled.

The two disjunctions shown in Equations \eqref{eq:sla1_disjunction} and \eqref{eq:sla2_disjunction} can be reformulated to equivalent MILP constraints using the associated binary variable $w$ as expressed in Equations \eqref{eq:sla1_reformulation} and \eqref{eq:sla2_reformulation} respectively.


\begin{subequations} \label{eq:sla_reformulations} 
\begin{align} 
{B}_{mst}^{SLA}w_{mst} &\leq {B}_{mst} \leq {B}^U_{mst}w_{mst} & \quad \forall \ s \in \mathcal{S}^{SLA}, \ m \in \mathcal{M}_s, \ t \in \mathcal{T} \label{eq:sla1_reformulation}\\ 
{B}_{mst}^{SLA}w_{mst} &\leq \sum_{t'=t}^{t+\tau^{SLA}_{mst}}{B}_{mst'} \leq {B}^U_{mst}w_{mst} & \quad \forall \ s \in \mathcal{S}^{SLA}, \ m \in \mathcal{M}_s, \ t \in \mathcal{T}: t+\tau^{SLA}_{mst} \leq |\mathcal{T}| \label{eq:sla2_reformulation} \end{align} \end{subequations}

This discussion only applies for the specific set of suppliers $\mathcal{S}^{SLA}$ that require service level agreements where $\mathcal{S}^{SLA} \subseteq \mathcal{S}$. 

The model that integrates SLA includes Equations \eqref{eq:delay}-\eqref{eq:final_time}, \eqref{eq:unmet_mip}-\eqref{eq:bounds}, and either Equation \eqref{eq:sla1_reformulation} or \eqref{eq:sla2_reformulation}.
As above, this extension requires additional binary variables and big-M constraints, which may negatively affect the solution time. These effects will be studied in Section \ref{sec:scalability}.

\subsubsection{Final Inventory Deviation (FID)}\label{sec:final_inventory} 

As mentioned earlier, providing a notion of recovery from the disruption is key 
to resume ``normal'' operation.
In the base model, this notion was given by Equation \eqref{eq:final_time} as a hard terminal constraint that matched pre and post-disruption states exactly.
These constraints can sometimes be too stringent and may potentially yield an unnecessary amount of late deliveries and cancellations. 
The hard constraints can be relaxed by replacing them with an L1 penalty in the objective to the final inventory deviation (FID) as:
\begin{align}
    Deviation_{mn} = |I_{mn}^0 - I_{mn|\mathcal{T}|}| \quad \forall \ n \in \mathcal{P} \cup \mathcal{W}, \ m \in \mathcal{M}_n  \label{eq:deviation_abs}
\end{align}
which can be reformulated as:
 \begin{subequations}
 \label{eq:deviation_refs}
 \begin{equation}
    Deviation_{mn} \geq I_{mn}^0 - I_{mn|\mathcal{T}|} \quad \forall \ n \in \mathcal{P} \cup \mathcal{W}, \ m \in \mathcal{M}_n  \label{eq:deviation_ref1}
\end{equation}
\begin{equation}
      Deviation_{mn} \geq I_{mn|\mathcal{T}|} - I_{mn}^0 \quad \forall \ n \in \mathcal{P} \cup \mathcal{W}, \ m \in \mathcal{M}_n  \label{eq:deviation_ref2}
\end{equation}
\end{subequations}
and which in turn, can replace Equation \eqref{eq:final_time} in the formulation.
Furthermore, the total deviation is calculated as:
\begin{equation}
      TotalDeviation = \sum_{n \in \mathcal{P} \cup \mathcal{W}} \sum_{m \in \mathcal{M}_n} \lambda^{Dev}_{mn} Deviation_{mn} \label{eq:deviation_obj}
\end{equation}
where a penalty ($\lambda^{Dev}$) incentivizes the difference to be as small as possible, and it can be subtracted from the objective function.
Intuitively, a large penalty can be applied when prioritizing rapid recovery from the disruption.
Therefore, the formulation extended to account for a final inventory deviation penalty includes Equations \eqref{eq:delay}-\eqref{eq:fixing_to_zero}, \eqref{eq:unmet_mip}-\eqref{eq:bounds}, and \eqref{eq:deviation_refs}-\eqref{eq:deviation_obj}.

The penalty cost is indexed by location and material enabling us to prioritize which products and places more urgently require a base inventory to resume their normal operation. 
It is worth noting that while this approach can help achieve solutions with fewer backorders and cancellations, certain assumptions (related to a gradual shrinkage of the deviations) must hold to ensure stability in a rolling horizon framework. However, addressing these assumptions is outside the scope of this work, therefore we refer the readers to Rawlings \citep{rawlings2017model} for more information on this subject.

\subsubsection{Negative Inventory Deviation (NID)} \label{sec:negative_dev}
The base model considers hard lower bounds on inventory variables. In this section, we extend the model with soft bounds to allow for more flexible responses. 
Contrary to Equation \eqref{eq:deviation_abs} where a symmetric L1 penalty was enforced as constraint, here we only impose a penalty when the inventory level is below a threshold given by the lower bound ${I^L}$. 
We refer to this one-sided deviation as negative inventory deviation (NID).
This alternative allows the model to go below the determined threshold at the cost of a penalty ($\lambda^{K}$). Naturally, larger values for $\lambda^{K}$ model a more conservative approach were inventory values above the threshold value are preferred. 
The logic is captured by the disjunction: 
\begin{equation}
    \left[\begin{gathered}
        Z_{mnt}\\
       I_{mnt} - {I}^L_{mnt} \leq 0 \\
       K_{mnt} = I_{mnt} - {I}^L_{mnt}
    \end{gathered}\right] \vee
    \left[\begin{gathered}
        \neg Z_{mnt}\\
       I_{mnt} - {I}^L_{mnt} \geq 0 \\
       K_{mnt} = 0
    \end{gathered}\right] \quad \forall \ n \in \mathcal{P}\cup\mathcal{W}, \ m \in \mathcal{M}_n, \ t \in \mathcal{T}
    \label{eq:invdev_disjunction}
\end{equation}
where $Z$ is the Boolean variable that indicates if there is a violation below the threshold, and $K$ represents the total negative deviation.
We also use here an MILP reformulation for the disjunction. 
First, we note that the original upper bounds on inventory posed in Equation \eqref{eq:I_bounds} are still required. Therefore, we can derive the following bound for the deviation:
\begin{equation}
       - {I}^L_{mnt} \leq I_{mnt} - {I}^L_{mnt} \leq {I}^U_{mnt} - {I}^L_{mnt}  \quad \forall \ n \in \mathcal{P}\cup\mathcal{W}, \ m \in \mathcal{M}_n, \ t \in \mathcal{T} \label{eq:inv_bound_k_bounds}
\end{equation}

\begin{figure}[!htb]
	\includegraphics[width=0.65\textwidth]{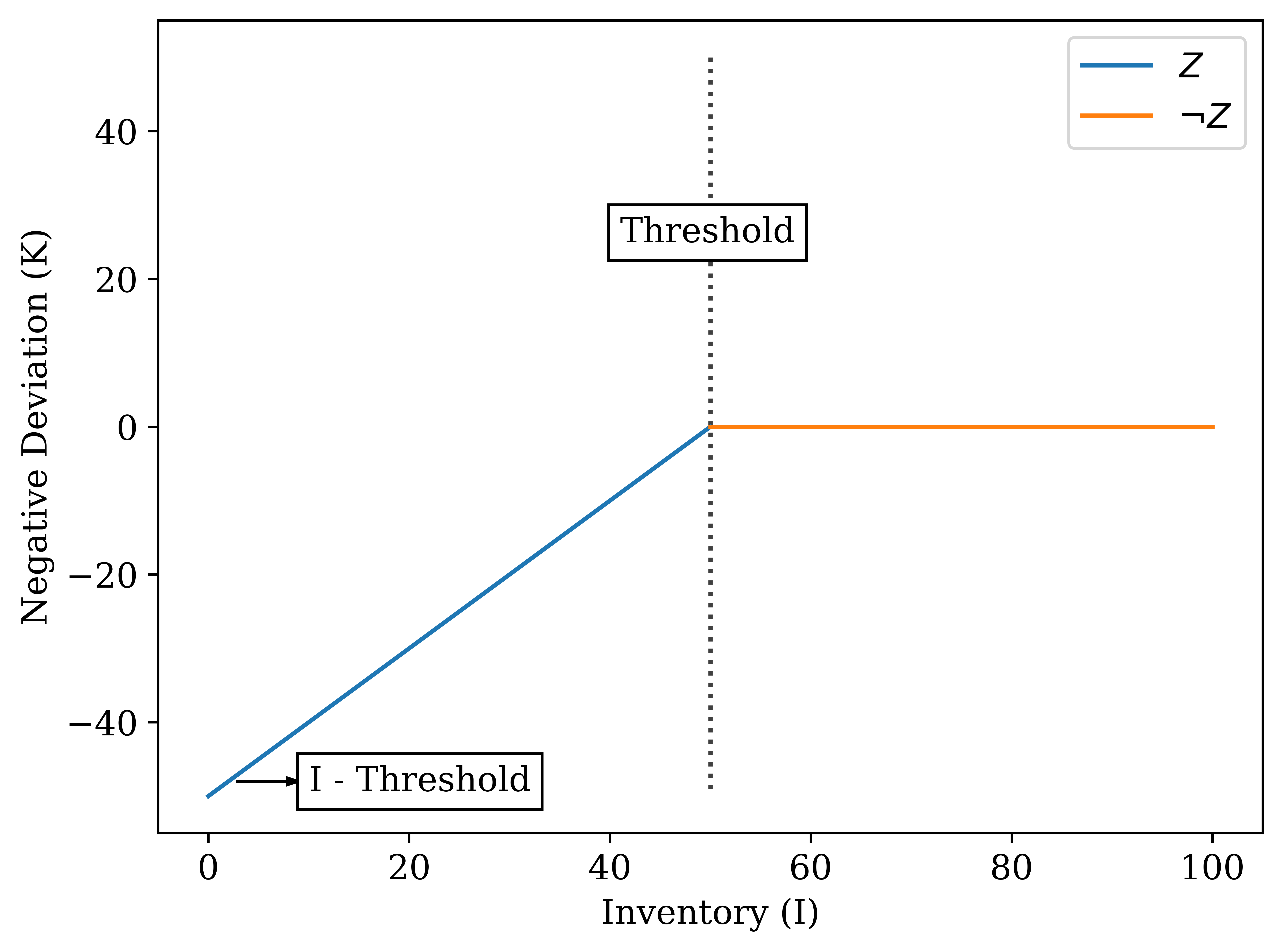}
	\centering
	\caption{Negative deviation from the threshold ${I}^L$ as a function of the inventory $I$ of material A in Plant 1. The inventory is bounded as $I \in [0,100]$ and a constant threshold  ${I}^L=50$ is imposed.}
	\label{fig:relu}
\end{figure}

The disjunction indicates that the negative deviation ($K$) behaves like a rectified linear unit (ReLU) reflected over the $y$-axis and displaced by ${I}^L$. Figure \ref{fig:relu} illustrates this penalty as a function of inventory. We reformulate the disjunction using the ReLU formulation from \citet{ceccon2022omlt}, reflected over the $y$-axis and displaced by the threshold, as follows

\begin{subequations}
\label{eq:relus}
\begin{equation}
     K_{mnt} \leq I_{mnt} - {I}^L_{mnt}  \quad \forall \ n \in \mathcal{P}\cup\mathcal{W}, \ m \in \mathcal{M}_n, \ t \in \mathcal{T}  \label{eq:relu1}
\end{equation}
\begin{equation}
    K_{mnt} \geq I_{mnt} - {I}^L_{mnt} + ({I}^L_{mnt}-{I}^U_{mnt})(1-z_{mnt})  \quad \forall \ n \in \mathcal{P}\cup\mathcal{W}, \ m \in \mathcal{M}_n, \ t \in \mathcal{T}  \label{eq:relu2}
\end{equation}
\begin{equation}
    K_{mnt} \geq -{I}^L_{mnt}z_{mnt}  \quad \forall \ n \in \mathcal{P}\cup\mathcal{W}, \ m \in \mathcal{M}_n, \ t \in \mathcal{T}  \label{eq:relu3}
\end{equation}
\begin{equation}
    K_{mnt} \leq 0  \quad \forall \ n \in \mathcal{P}\cup\mathcal{W}, \ m \in \mathcal{M}_n, \ t \in \mathcal{T}  \label{eq:relu4}
\end{equation}
\end{subequations}

The deviation along the time horizon is aggregated and penalized to obtain the total cost of inventory bounds violation:
\begin{align}
      TotalNegativeDeviation = \sum_{n \in \mathcal{P} \cup \mathcal{W}} \sum_{m \in \mathcal{M}_n} \sum_{t \in \mathcal{T}} \lambda^{K}_{mnt} K_{mnt} \label{eq:underpass_obj}
\end{align}
which has to be added to the objective of the base model to penalize lower bound violations of the threshold. Therefore, the model incorporates a lower bound penalization for the inventory is given by Equations
\eqref{eq:delay}-\eqref{eq:final_time}, \eqref{eq:unmet_mip}-\eqref{eq:bounds}, and \eqref{eq:inv_bound_k_bounds}-\eqref{eq:underpass_obj}. 

In the section above, we established the base model as a powerful formulation for reactive responses in supply chain and manufacturing networks. 
Additionally, we expanded the model to accommodate scenarios with relaxed initial assumptions, yielding more adaptable responses. 
In Section \ref{sec:scalability} we will discuss the implications of these extensions on the computational performance.
In the following section, we will demonstrate how to integrate various disruptions that occur in supply chains into the model.

\section{Disruption Modeling} \label{sec:disruption_modeling}

The proposed framework enables the modeler to incorporate various types of disruptions that may arise in a supply chain system. 
Generally, a disruption denotes a deviation of the system from its nominal operation.
This study classifies disruptions into two primary categories: immediate disruptions and scheduled disruptions. 
Immediate disruptions occur unexpectedly or in an unplanned manner, including events like equipment failures, roadblocks, or supply shortages. 
The model formulation is designed to generate a reactive schedule in response to these disruptions, allowing the system to adapt to the unexpected changes.
Conversely, scheduled disruptions are events for which abnormal circumstances are known in advance. 
This category includes instances such as regular maintenance or planned shutdowns. 
In these scenarios, the formulation can offer proactive policies, enabling the network to prepare for the disruption ahead of time and operate effectively throughout the disruption period.


Characterizing a disruption involves two parts. First, identifying the disrupted horizon meaning when the disruption starts and ends. The start can be at $t=0$ for an unplanned disruption. Similarly, the end can be at period $t=|\mathcal{T}|$ if the end of the  disruption is unclear or if the system is known to be permanently disrupted.
Second, the set of disrupted parameters must be identified, which captures the negative impacts on the system, causing deviations from normal operation throughout the disrupted horizon.

A significant portion of disruptions occurring in real-life systems can be mapped, in some manner, to a reduction of available capacity. 
Examples of this concept include a roadblock impacting edge capacity, a supply shortage reducing purchasing capacity, or an equipment failure decreasing production capacity. 
By indexing both upper and lower bounds of continuous variables over time, bound profiles can be generated to model disruptions by tightening the bounds.

\begin{figure}[!htb]
\centering
\begin{subfigure}{0.49\textwidth}
    \includegraphics[width=1.1\textwidth]{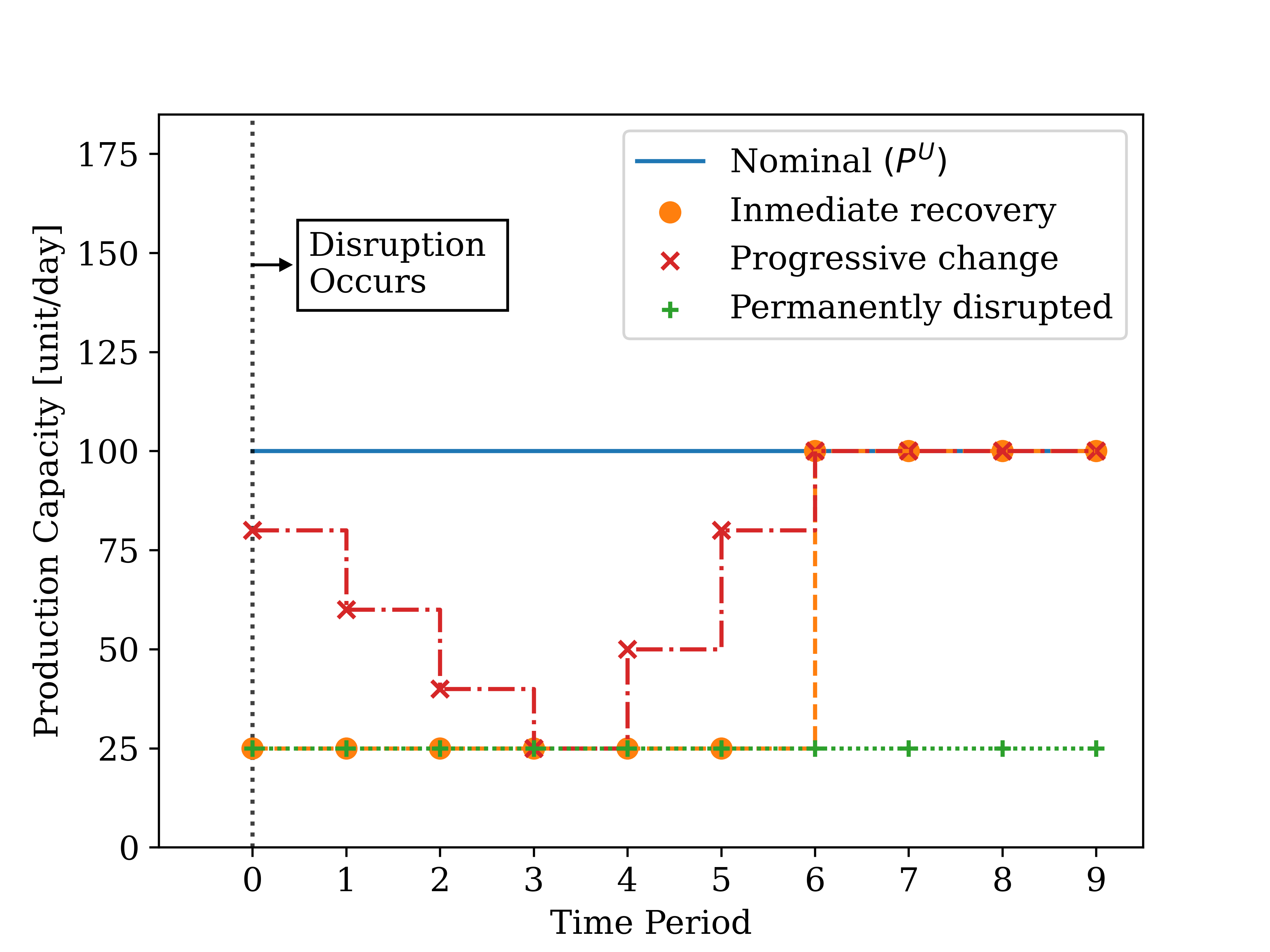}
    \caption{Bound profiles to model unexpected (immediate) disruptions.}
    \label{fig:variable_bound_profiles_1}
\end{subfigure}
\hfill
\begin{subfigure}{0.49\textwidth}
    \includegraphics[width=1.1\textwidth]{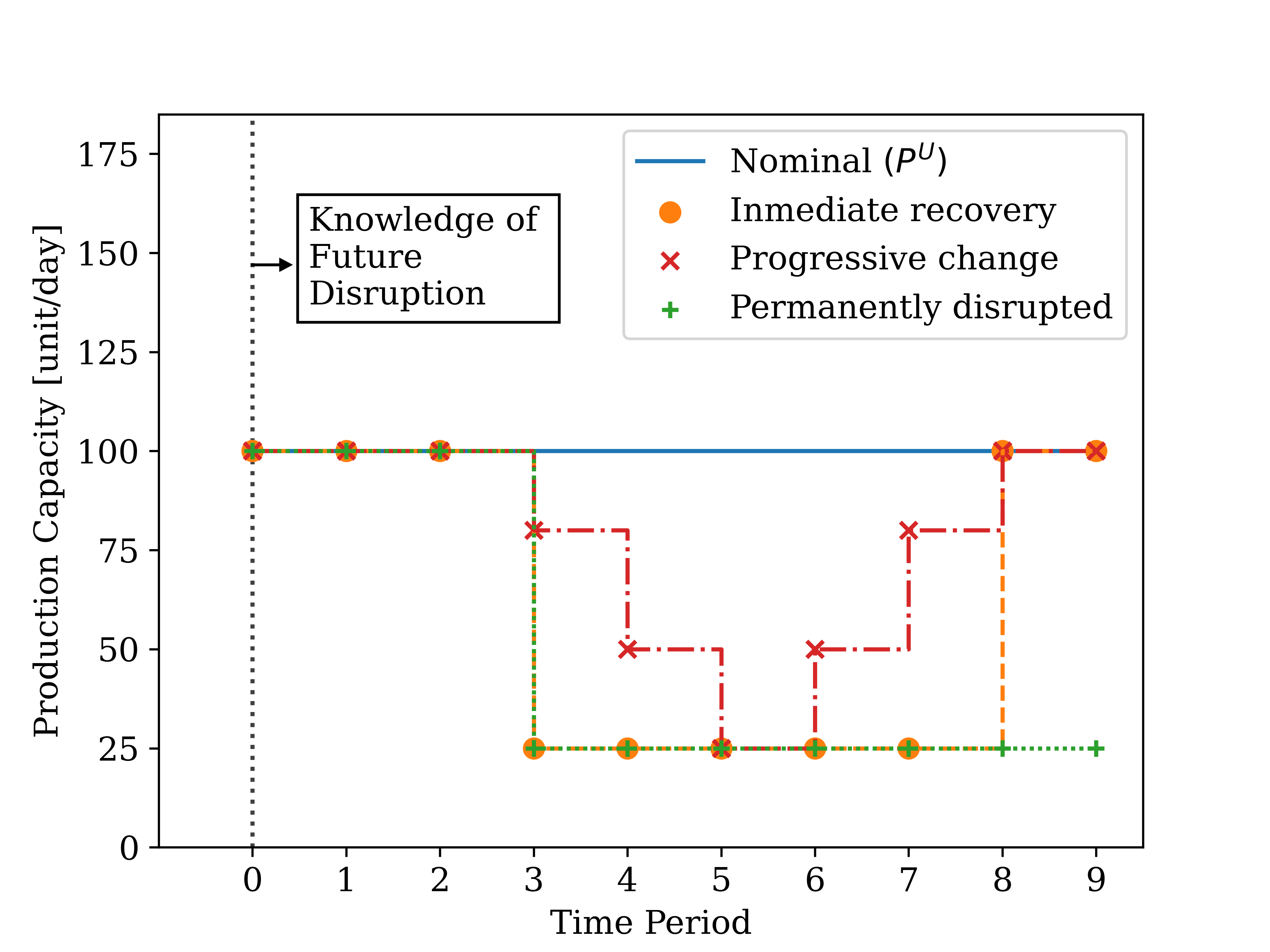}
    \caption{Bound profiles to model known (scheduled) disruptions.}
    \label{fig:variable_bound_profiles_2}
\end{subfigure}
        
\caption{Potential production upper bound profiles for the production of a recipe  to model a disruption that reduces the capacity by 75\% in both immediate and scheduled cases.}
\label{fig:variable_bound_profiles}
\end{figure}

Consider the example where there is a reduction in plant production, hence the disrupted parameter corresponds to production bounds.
Figure \ref{fig:variable_bound_profiles} illustrates how both immediate and scheduled disruptions can be modeled by designing a bound time profile.
These profiles not only encode information about the time of disruption and recovery, but also about how the capacity precisely decays and increases over time.

The nominal upper bound, depicted in blue, represents the capacity of the current operation, which serves as a baseline for the undisrupted system. 
In certain disruptions, the system recovers immediately after the event has passed, as shown in orange, where no additional time to return to its usual capacity.
In other cases, disruptions may progressively tighten and then relax the bound over time, illustrated in red, offering the modeler flexibility to characterize the disruption dynamically. 
In green are modeled scenarios where the system is assumed to be permanently disrupted, either due to a severe disruption, or a lack of knowledge about when the disrupted horizon will end.
Furthermore, both bounds can collapse to a single value at a specific time if necessary. 
For example, having the upper bound go to zero can model a complete shutdown in plant production.

Other disruptions can be modeled by changing parameters other than variable bounds. For instance, if shipping through a route takes longer due to external factors, the reactive response might involve finding alternative routes, and adjusting downstream operations to handle the lead time disruption. 
In our approach, transportation lead times can vary over time.
To model an increase in transportation delay, simply adjust the nominal lead times.
This discussion can also be extended to production delays as proposed in Section \ref{sec:production_time}.


In supply chain operations, there are instances where economic values, such as costs and prices change often due to external factors that affect the various stakeholders involved in the operation. 
Examples include market fluctuations impacting product prices, or supply shortages affecting not only supplier capacity but also the cost of the product. 
To model these changes in economic parameters, it suffices to impose the disrupted prices and costs.
By doing so and obtaining an optimal reactive schedule, decision-makers can assess how significantly the optimal reaction would differ from the usual operation.


There is one case of disruption we consider in this work that aligns with what most of the literature has focused on when it comes to supply chain networks under uncertainty: changes in the demand.
Unlike most literature that provides fixed policies to address sudden changes or spikes in demand, our reactive approach aims to optimize the operation of the system in response to unexpected orders. 
If a new order is placed, our approach evaluates how to readjust the system to integrate the incoming order.
Here the model may decide to deliver the order late or cancel it to avoid inducing stress on the system. 
In this case, the set of disrupted parameters must not only include the size of the order but also the associated late delivery cost profile and cancellation cost of the order place at period $t'$ for the specific material $m'$ and customer $c'$.
Furthermore, if there are specific reasons (such as customer importance or contractual obligations) for which the new order cannot be delivered late or canceled, this policy can be enforced by:
 \begin{subequations}
 \label{eq:new_order}
 \begin{equation}
    U_{m'c't} = 0 \quad \forall \ t \in \mathcal{T} \label{eq:new_order_late}
\end{equation}
\begin{equation}
      y_{m'c't'} = 0 \label{eq:new_order_cancel}
\end{equation}
\end{subequations}
in which case the final state of model has to be extended, as proposed in Section \ref{sec:final_inventory} to guarantee feasibility given that an order is being forced to be delivered.

This proposed approach is capable of integrating various types of disruptions that occur in real life in an equation-oriented manner to provide an optimal response. 
Impacts in operations are modeled independently by design, such that the modeler determines if multiple disruptions are present in a single event.
Furthermore, the proposed method gives a way to bridge most of the existing literature to account for consequences of unexpected demands. In the following sections of the paper, we  present a case study and evaluate the computational performance of the model.


\section{Motivating Case Study} \label{sec:case_study}
\begin{figure}[!htb]
	\includegraphics[width=0.90\textwidth]{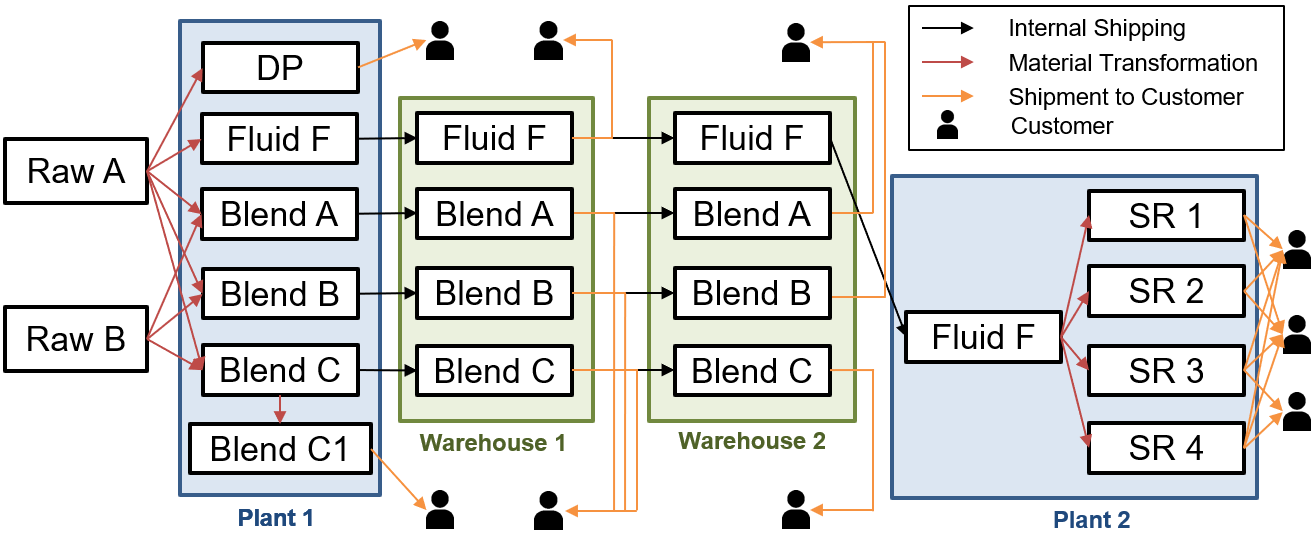}
	\centering
	\caption{ Supply chain and manufacturing network topology for the motivating case study}
	\label{fig:network}
\end{figure}

The supply chain network we consider is motivated by a part of the silicone rubber business of The Dow Chemical Company. 
It manufactures ten different products, and consists of two raw material suppliers, two plants, two warehouses, and nine customers that are satisfied from both manufacturers and warehouses. 
We consider the operation over a four month period ($|T| = 120$ days) and the structure of the motivating example is presented in Figure \ref{fig:network}.

The first plant (Plant 1) uses the two different raw materials (Raw A and Raw B) to make six products (DP, Fluid F, Blend A, Blend B, Blend C and Blend C1). 
In this plant, two products (labeled in the figure as DP and Blend C1) are sold to customers directly, while the rest of the products are packed in steel drums and shipped to Warehouse 1 through road transportation. 
Note that Blend C is both an intermediate material in the production of Blend C1 as well as a final product that is delivered to Warehouse 1. 
Once in Warehouse 1, these four products can be sold directly to local customers or shipped to Warehouse 2 through marine or air transportation. 
There, three products (labeled in the figure as Blend A, Blend B, and Blend C) are sold directly to local customers. 
Fluid F is shipped to Plant 2 to be further processed into four different final products (SR 1 - SR 4) and sold to local customers.
All nine customers have independent demand frequencies and order quantities.

Due to its size and complexity, the network can be exposed to several unplanned events that cause disruptions in its operation. This particular network can face, for example, four different types of disruptions that could negatively impact the overall operation.
\begin{itemize}
    \item Raw material shortage: This can occur when an upstream supplier is unable to deliver the amount of raw material originally requested by the company for scheduled operation. 
    \item Plant equipment failure: This causes a reduction or complete loss of production at a plant facility, impacting downstream consumers and inventory levels for both upstream and downstream elements.
    \item Route unavailability: A given arc in the supply chain may not capable of transporting any material for a specific time window. In this case, the network becomes partially disconnected and alternative routing options may be needed.
    \item Logistics resource limitation: The transport of material between supply chain elements can be impacted by availability of appropriate logistics resources (e.g., packing materials like steel drums or availability of drivers for truck transport). Furthermore, a limitation in the resources used for storage can affect the capacity in warehouses.
\end{itemize}
These examples can cause a reduction in available material for downstream processes and delays in deliveries of customer orders. However, these can also impact upstream operations, and cause inventory management challenges. For example, a plant failure causes an increase in raw material inventory that can impact upstream suppliers and producers. Hence, an integrated approach that considers the network holistically is desired.


\section{Results} \label{sec:results}

In the following subsections, we present the computational results of the proposed model. 
We begin by presenting the optimal reactive schedules generated by the model for two scenarios.
This provides insight into how the model determines an integrated schedule for supply chain operations and demonstrates how disruptions can propagate beyond the directly affected nodes.
Then, we investigate how the computational solution time is influenced by increasing the size of the model, particularly by extending the time horizon and examining the impact of different model extensions.  
This section provides details on the size of the models studied, including the number of binary variables, continuous variables, constraints, and non-zero elements.
Subsequently, we evaluate the different optimal operations obtained with the model across various characterizations of the same disruption.
Finally, we conduct a sensitivity analysis on the cancellation and late delivery penalties, which are key parameters that impact optimal decisions.

The proposed solution approach, which includes all the extensions, was implemented in \texttt{Pyomo} \citep{bynum2021pyomo} and the resulting MILP was solved with \texttt{Gurobi v9.5.1} as the solver on a Linux machine with 8 Intel\textregistered \, Xeon\textregistered \, Gold 6234 CPUs running at 3.30 GHz with 8 total hardware threads and 1 TB of RAM running Ubuntu. 
All the data presented from the Dow Chemical Company has been anonymized, therefore, the values and magnitudes showcased in this section do not reflect the actual operational figures of the company. 
The source code and sample data for illustration are freely available in the following open-source \hyperlink{https://github.com/dovallev/supply_chain_disruptions}{Github repository}\footnote{https://github.com/dovallev/supply\_chain\_disruptions}.

\subsection{Optimal Reactive Operation} 

In this section, we demonstrate the proposed formulation on a single case study and
 analyze the optimal solution.
 We consider the supply chain and manufacturing network outlined in Section \ref{sec:case_study}, and we compare the solution from a single disruption scenario with the solution in the undisrupted case. 
 The resulting models have 17,345 continuous variables, 497 binary variables, and 15,480  constraints for the base model with a solution time of 0.23 seconds.
 The disruption involves the unplanned failure (meaning that it occurs at time zero) of two out of three reactors at Plant 1, which produce DP, Fluid F, and Blend A. This reduces their production capacity to one-third of the nominal value for two months. We examine the optimal response of the base model from Section \ref{sec:base_model} over the next four months (120-day discretization), focusing on the transformation of Raw A into Fluid F and its subsequent use in producing SR 3. While our discussion centers on these specific materials, it is important to note that the model provides a comprehensive schedule for all materials in the network, with the disruption impacting the entire operation due to the critical nature of Raw A.

\begin{figure}[!htb]
\centering
\begin{subfigure}{0.49\textwidth}
    \includegraphics[width=1\textwidth]{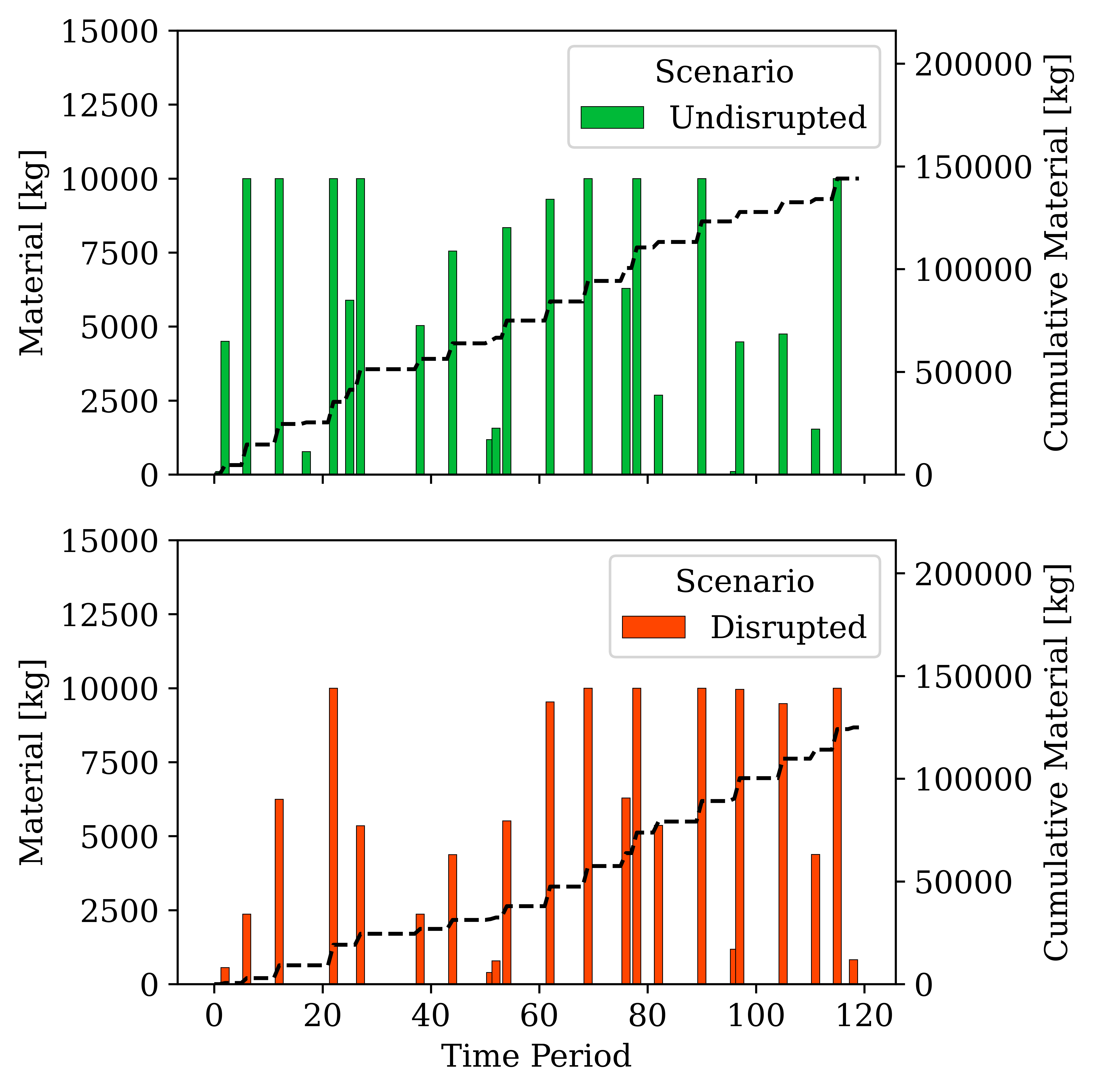}
    \caption{Total procurement of Raw A from supplier.}
    \label{fig:toy_rawa_buy}
\end{subfigure}
\hfill
\begin{subfigure}{0.49\textwidth}
    \includegraphics[width=1\textwidth]{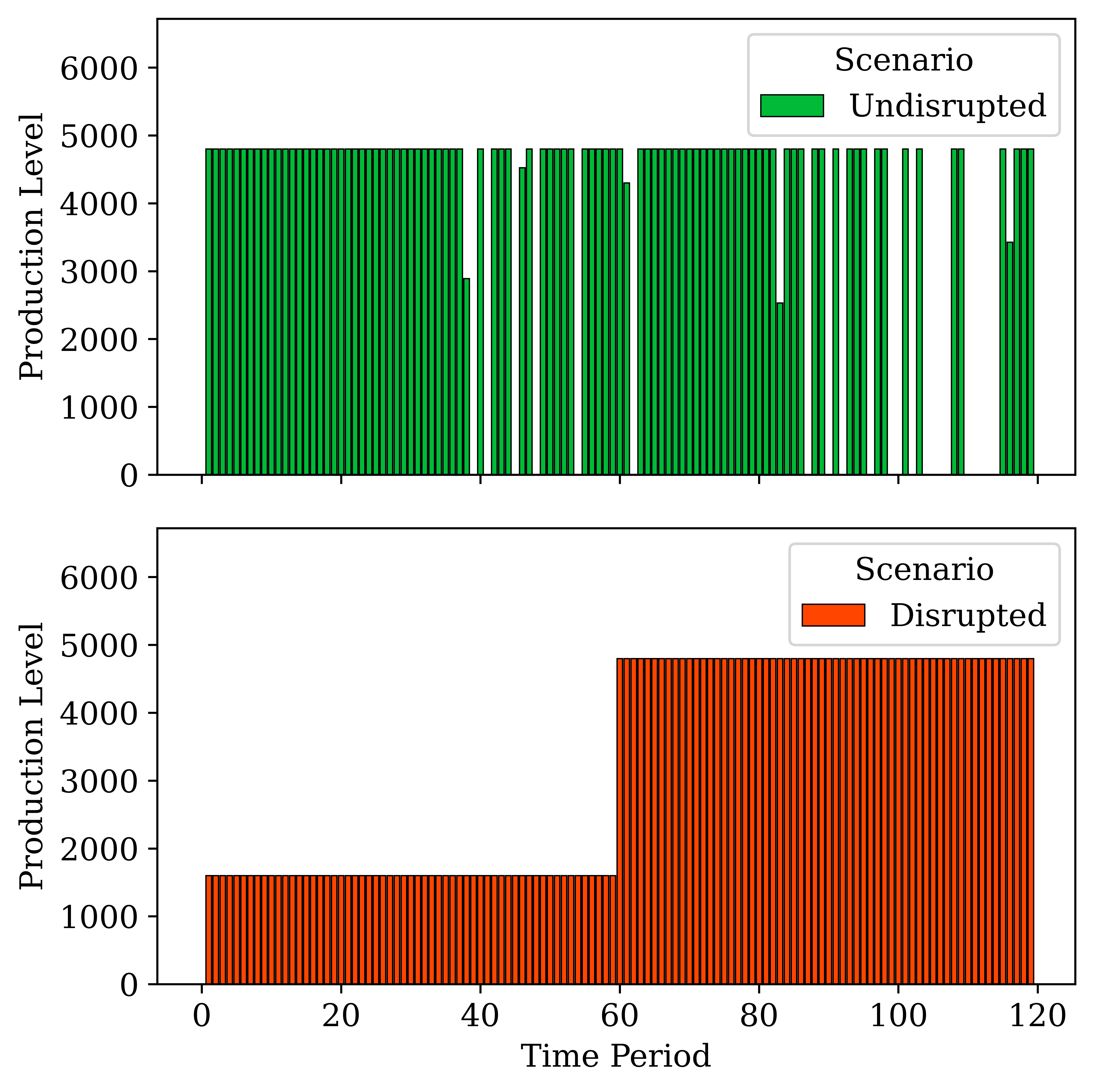}
    \caption{Production of Fluid F from Raw A in Plant 1.}
    \label{fig:toy_fluif_prod}
\end{subfigure}

\caption{Procurement and production levels in Plant 1 for Raw A and Fluid F.}
\label{fig:toy_rawa_to_fluidf}
\end{figure}

 Figure \ref{fig:toy_rawa_to_fluidf} illustrates the transformation of Raw A into Fluid F that takes place in Plant 1 for both undisrupted and disrupted scenarios. More specifically, Figure \ref{fig:toy_rawa_buy} presents the procurement schedule of Raw A from an external supplier, 
while Figure \ref{fig:toy_fluif_prod} depicts the production schedule for Fluid F from Raw A in Plant 1.
The disruption shifts larger procurement orders to the second half of the operation in contrast to the undisrupted case where most orders are placed in the first half. 
This can be seen best by comparing the slope of cumulative procurement line in each scenario: the undisrupted scenario shows faster procurement in the first half and the disrupted scenario shows faster procurement in second half.
This occurs because, with limited production capacity in the first half of the horizon, the model minimizes raw material procurement to avoid excess inventory. As production ramps up in the second half, increased raw material is then procured to meet the increased demand.
This shows how disruptions can impact upstream operations: a disruption at Plant 1 affects the production rate and subsequently affects the need for raw materials and the flow from the upstream supplier.

In the disrupted scenario, the production capacity for Fluid F is reduced to one-third, requiring the plant to operate at full capacity for the entire period to mitigate the disruption. In contrast, the undisrupted case allows for periods where Fluid F production is not at full capacity. 
Naturally, the plant operates at full capacity once production is restored, aiming to recover quickly and minimize the impact of the disruption.

\begin{figure}[!htb]
\centering
\begin{subfigure}{0.49\textwidth}
    \includegraphics[width=1\textwidth]{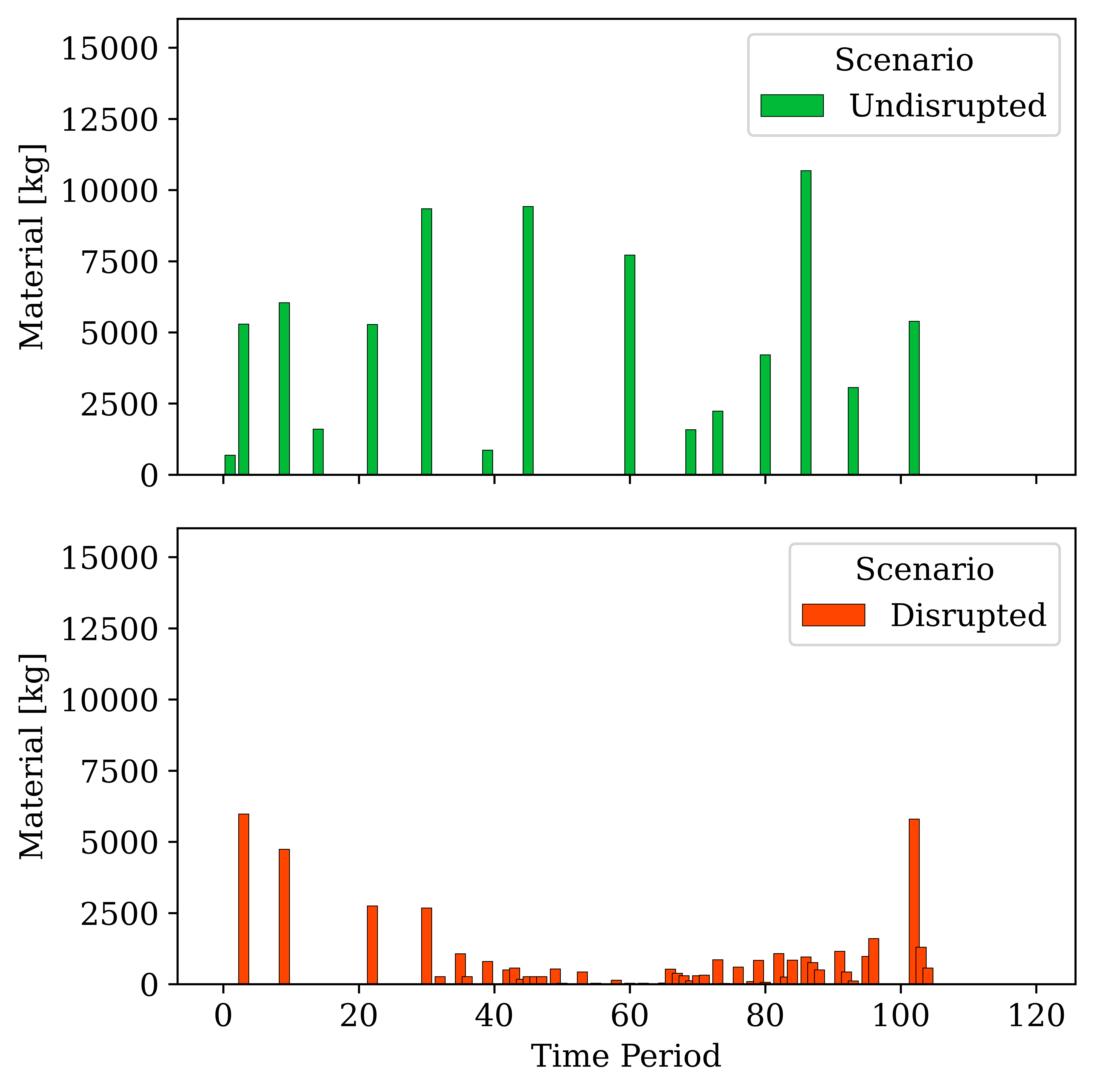}
    \caption{Shipment via truck.}
    \label{fig:toy_w1w2_truck}
\end{subfigure}
\hfill
\begin{subfigure}{0.49\textwidth}
    \includegraphics[width=1\textwidth]{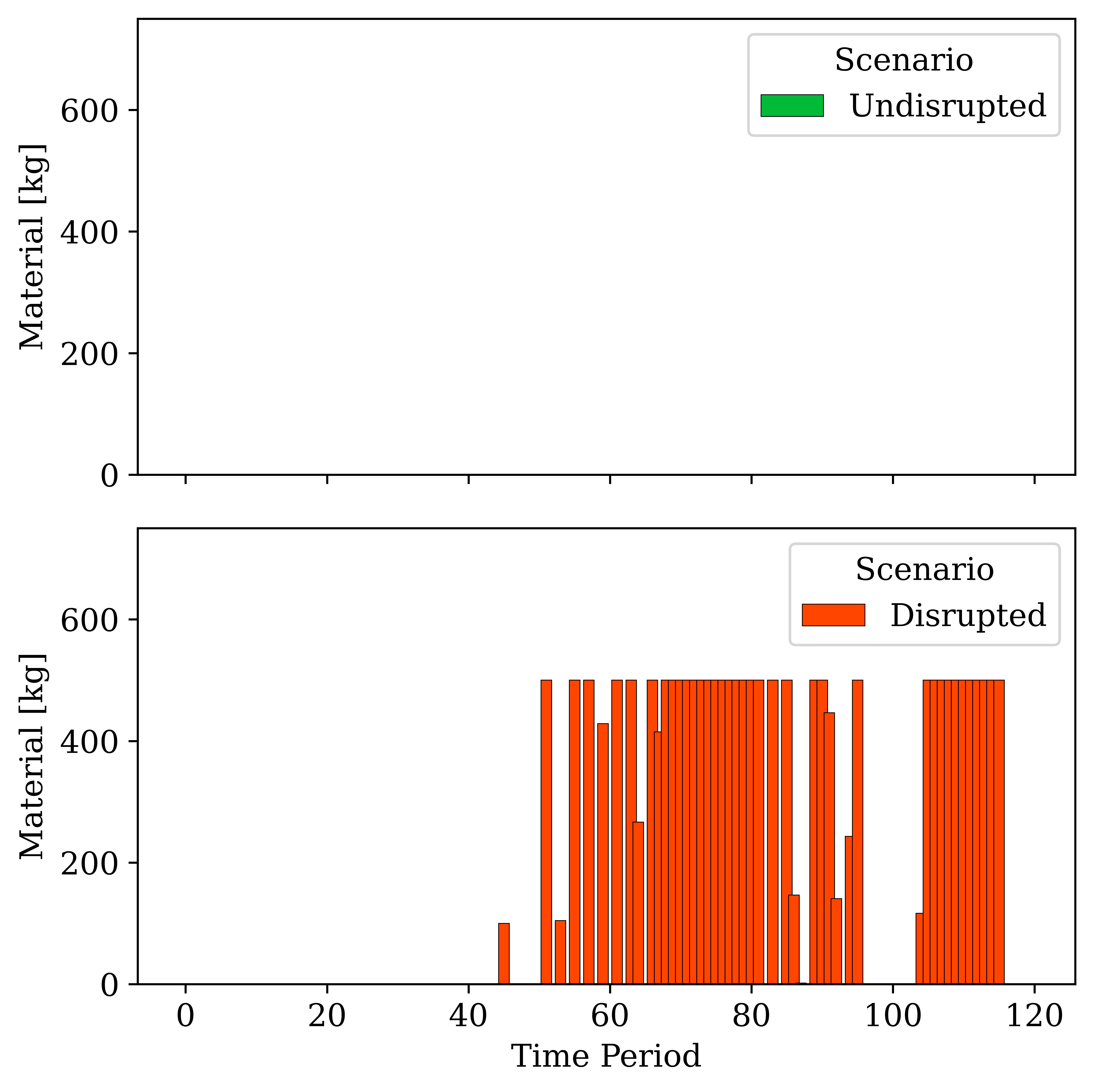}
    \caption{Shipment via air-freight.}
    \label{fig:toy_w1w2_air}
\end{subfigure}

\caption{Shipment of Fluid F from Warehouse 1 to Warehouse 2.}
\label{fig:toy_fluid_movement}
\end{figure}

Figure \ref{fig:toy_fluid_movement} illustrates how the disruption extends beyond the affected node. Specifically, Figures \ref{fig:toy_w1w2_truck} and \ref{fig:toy_w1w2_air} depict Fluid F shipments from Warehouse 1 to Warehouse 2 via truck and air freight, respectively. 
The disruption significantly impacts the optimal shipment of Fluid F. In the undisrupted case, large, sporadic shipments are sent via truck, while the disrupted case sees more frequent, smaller shipments. Notably, airfreight, unused in the undisrupted scenario, is employed from Warehouse 1 to Warehouse 2 in the second half of the disrupted case. 
This shift occurs as the operation tries to ``catch up'' using faster, but more expensive transportation modes.

\begin{figure}[!htb]
\centering
\begin{subfigure}{0.49\textwidth}
    \includegraphics[width=1\textwidth]{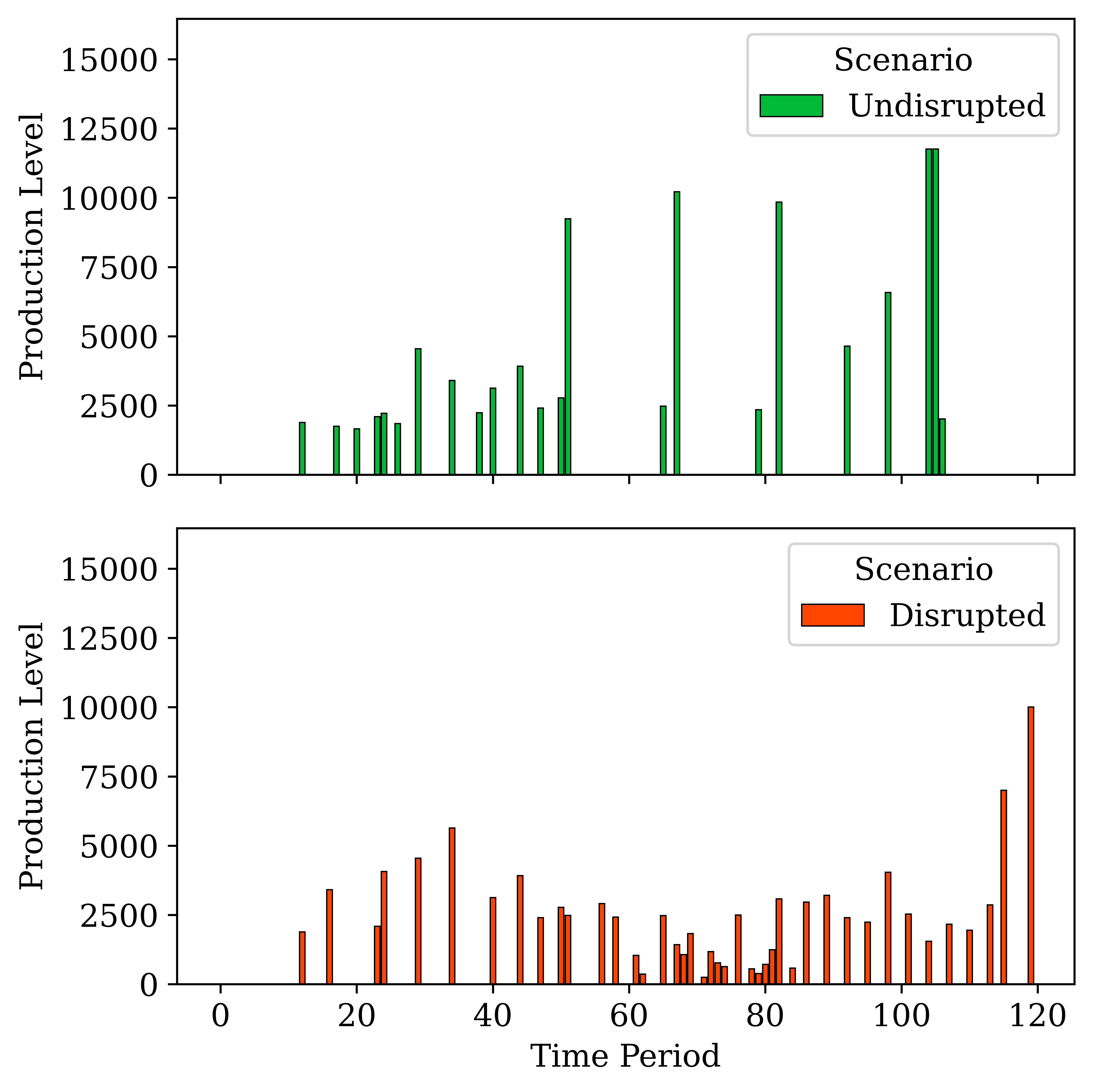}
    \caption{Production of SR 3 from Fluid F in Plant 2.}
    \label{fig:toy_sr3_prod}
\end{subfigure}
\hfill
\begin{subfigure}{0.49\textwidth}
    \includegraphics[width=1\textwidth]{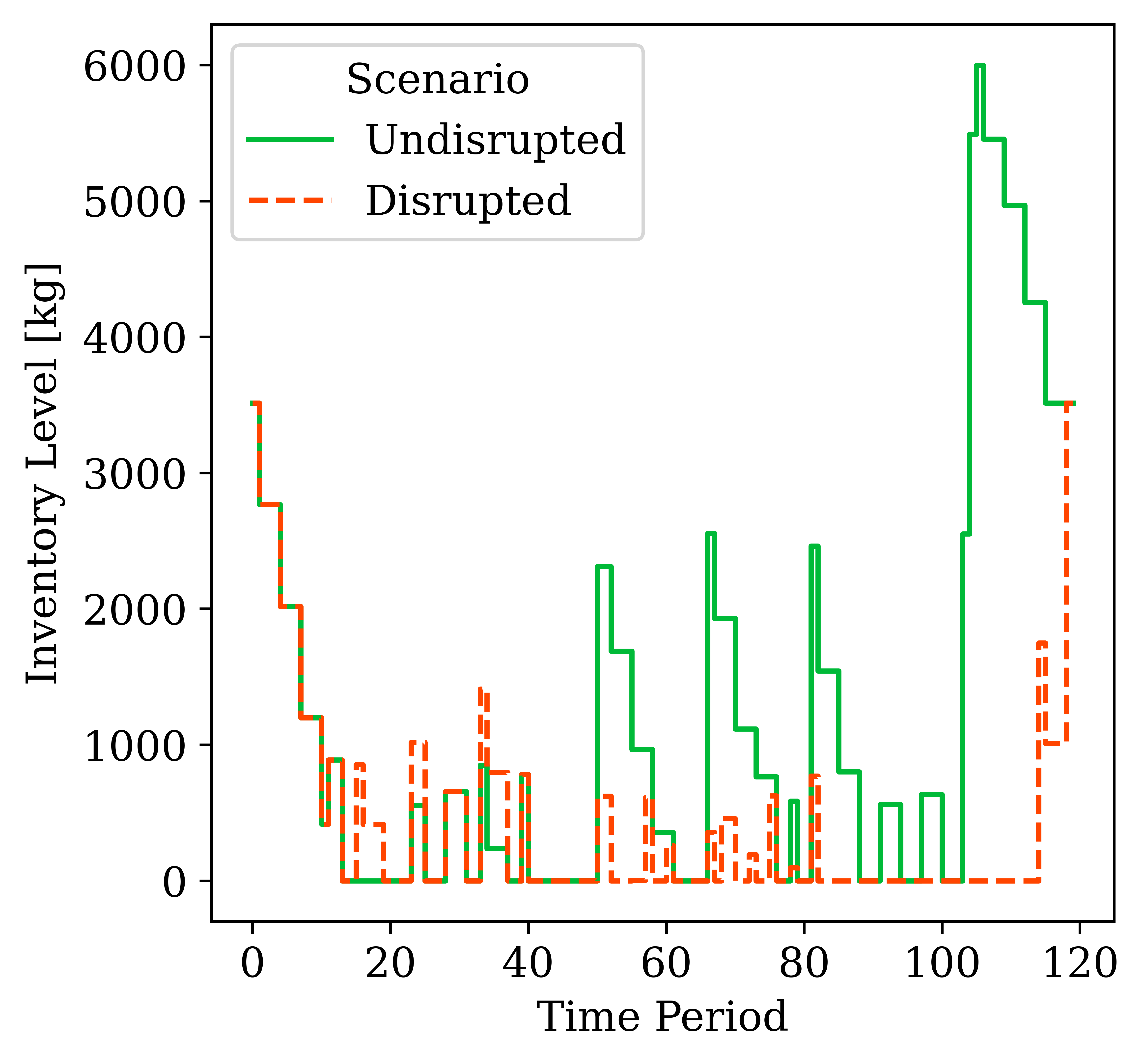}
    \caption{Inventory level of SR 3 in Plant 2.}
    \label{fig:toy_sr3_inv}
\end{subfigure}
        
\caption{Production and inventory levels of SR 3 in Plant 2.}
\label{fig:toy_sr3_prodinv}
\end{figure}

The reactor disruption not only directly impacts Fluid F production, but also affects other materials downstream. This phenomenon is illustrated by the production and inventory profiles of SR 3 at Plant 2, as shown in Figures \ref{fig:toy_sr3_prod} and \ref{fig:toy_sr3_inv}, respectively.
The SR 3 profiles reveal how the upstream disruptions propagate downstream. In the undisrupted scenario, production orders are larger and more sporadic, whereas, in the disrupted case, orders are smaller and more frequent. 
These behaviors directly reflect the  downstream propagation of the shipping patterns discussed in Figure \ref{fig:toy_w1w2_truck}. Since Fluid F is the raw material used in production, the production schedule in each scenario depends on the quantity and frequency of Fluid F shipments from upstream. 
Notably, the disrupted scenario demonstrates the ability of the model to generate an integrated and centralized response, as the strategy to recover quickly (using small, frequent shipments and production runs as they become available) is coordinated across the entire network.

The inventory profiles initially align in both scenarios, but they diverge over time as the disruption propagates through the network. 
This behavior occurs because the disruption takes time to propagate to this node in the network. This makes sense given that operations remain unaffected as they rely on existing inventory rather than upstream supplies. Notably, this approach in the initial stages is cost-effective, as utilizing existing stock is cheaper than shipping material from other nodes.
In the undisrupted case, inventory increases correspond to production periods with gradual depletion as external demand arises. In the disrupted scenario, we observe lower inventory levels and a higher turnover frequency.
In the latter case, the high turnover frequency results from an integrated operation that is trying to recover quickly and is unable to hold inventory for extended periods.

\begin{figure}[!htb]
\centering
\begin{subfigure}{0.49\textwidth}
    \includegraphics[width=1\textwidth]{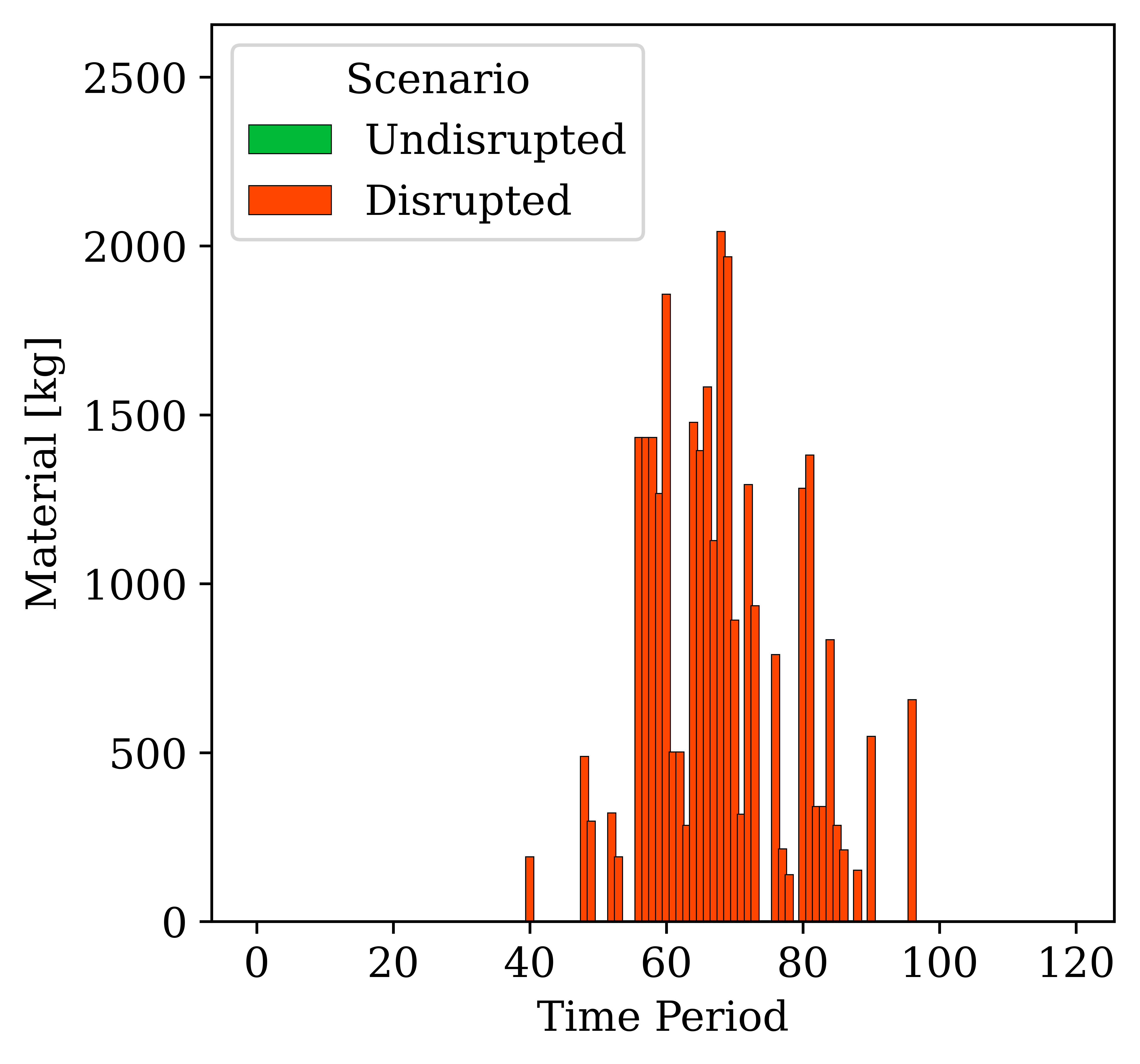}
    \caption{Total amount of material delayed.}
    \label{fig:toy_dem_unmet}
\end{subfigure}
\hfill
\begin{subfigure}{0.49\textwidth}
    \includegraphics[width=1\textwidth]{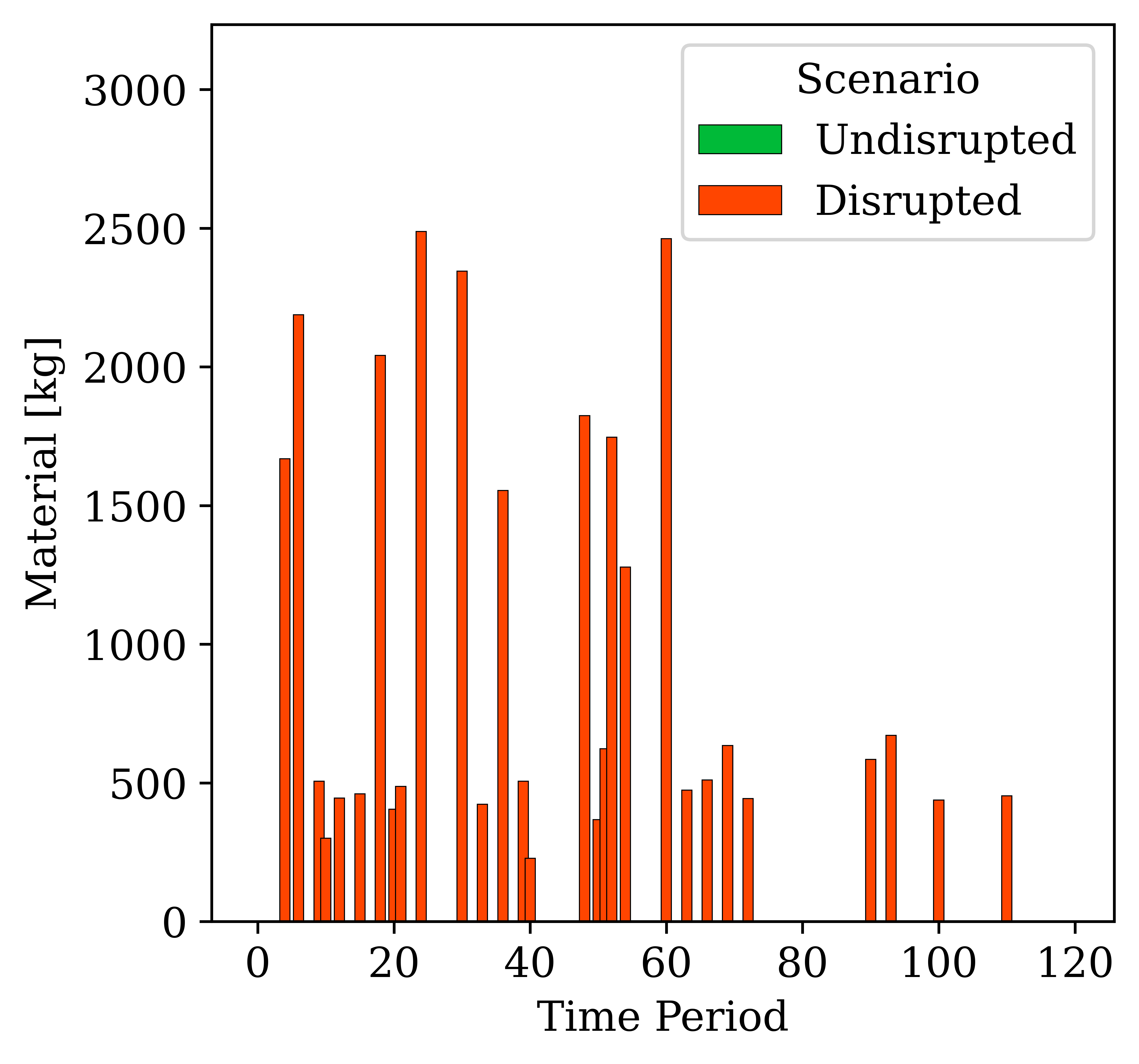}
    \caption{Total amount of material canceled.}
    \label{fig:toy_cancel}
\end{subfigure}
        
\caption{Total amount of material delayed and canceled throughout the operation.}
\label{fig:toy_unmetcancel}
\end{figure}
One of the most significant consequences of disruptions is their impact on demand coverage, affecting the timely delivery of orders to customers. Figure \ref{fig:toy_unmetcancel} illustrates the order management outcomes of the model. Specifically, Figure \ref{fig:toy_dem_unmet} shows the total amount of delayed material, while Figure \ref{fig:toy_cancel} displays the total amount of material from canceled orders. This analysis aggregates data across all 10 materials sold by the network.
Analyzing both figures, we observe that the undisrupted case involves no delays or cancellations, as expected for an operation designed for steady-state nominal conditions (hence, no green profiles are present in the figures).
In the disrupted scenario, the model responds by canceling large orders in the first half of the time horizon due to insufficient resources, avoiding the high costs of long delays. This early cancellation strategy is complemented by delaying orders in the middle of the time horizon, where there is less time to wait for resources, allowing the fulfillment of these orders with some delay. Notably, the overall order management response is integrated. These decisions are known immediately upon the disruption, enabling decision-makers to communicate with customers and negotiate alternative solutions in advance.

\subsection{Time Scalability} \label{sec:scalability}
\begin{figure}[!htb]
	\includegraphics[width=0.7\textwidth]{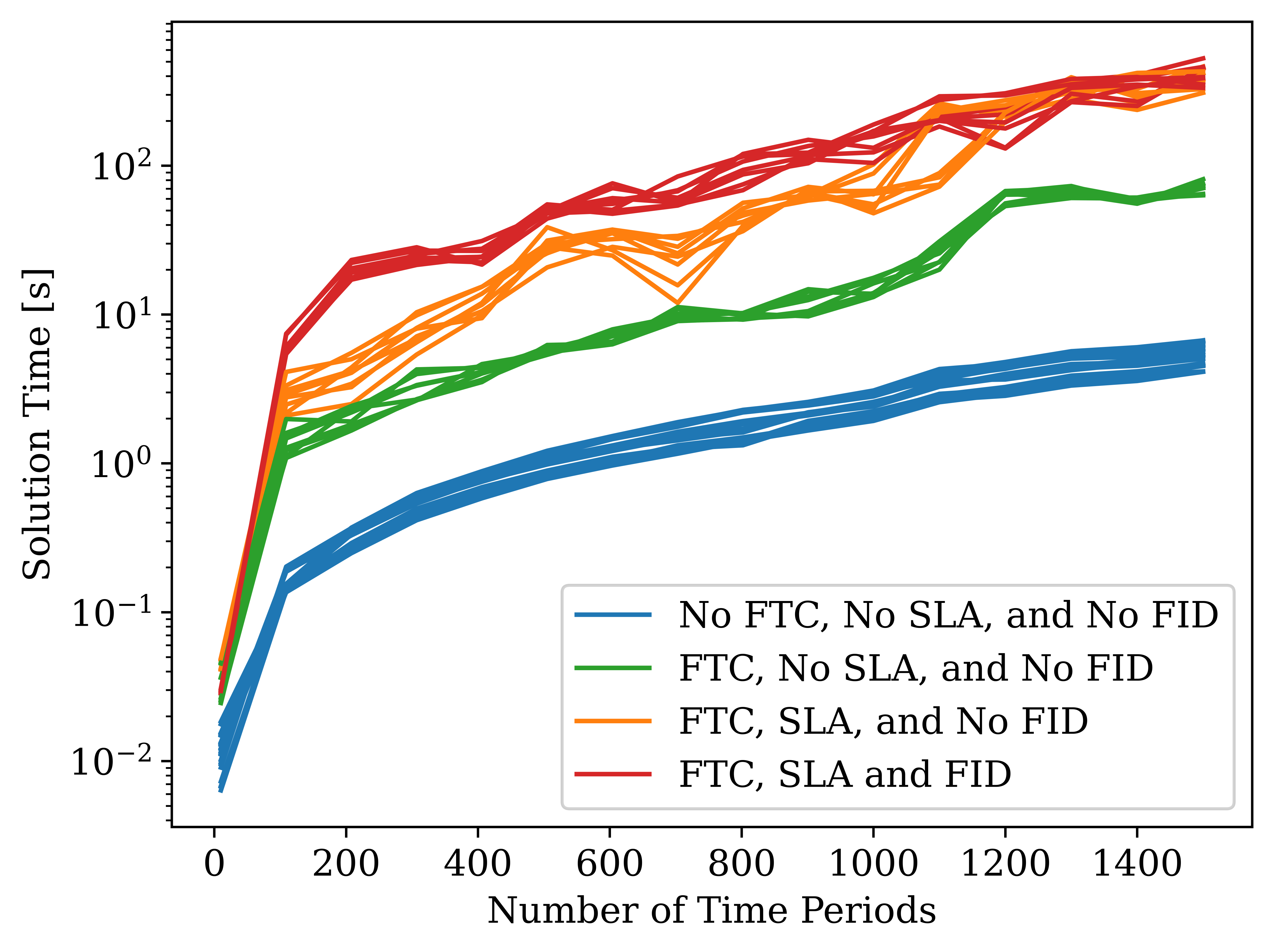}
	\centering
	\caption{Time scalability of the different model extensions.}
	\label{fig:time_scalability}
\end{figure}

Determining the number of time periods in the horizon is a modeling choice that can significantly impact the computational time. 
In practice, increasing the number of time periods can be done to increase the horizon or to increase the frequency of decisions.
In this section we study the impact that the number of time periods has on the solution time of the different extensions to an undisrupted model. 
We evaluated all possible combinations of models that can be created by including (or not) the extensions discussed in Section \ref{sec:extensions} yielding a total of 48 models. 
For each model type, a set of 50 replicas are constructed by varying order schedules. 
These were all solved for increasing number of time periods resulting in models with up to hundreds of thousands of variables and constraints as shown in Table \ref{tab:size_summary}. 
Figure \ref{fig:time_scalability} shows the average solution time of all the replicas where each individual line corresponds to a different model, summarizing over 290 hours of computing time.

\begin{table}[htbp]
\centering
\begin{tabular}{|l|c|c|c|c|}
\hline
\textbf{}                   & \textbf{Min} & \textbf{Max} & \textbf{Mean} & \textbf{Median} \\ \hline
\textbf{Binary Vars.}             &     25         &      98,231        &      26,540.81         &    19,370             \\ \hline
\textbf{Continuous Vars.}         &       1,199       &      266,666        &       122,111.84        &  121,294               \\ \hline
\textbf{Constraints} &        898       &     377,889        &       132,613.65        &       122,516          \\ 
\hline
\textbf{Non-zeros} &       2,421       &       881,549       &       332,432.11        &         319,199         \\ \hline
\end{tabular}
\caption{Model size summary statistics.}
\label{tab:size_summary}
\end{table}


While  48 different models were considered, we have identified several key extensions that have the strongest impact on solution time. 
Figure \ref{fig:time_scalability} includes the results for all models; however, we group them according to these high-impact extensions.
The y-axis of the figure represents computational time on a logarithmic scale, thus unit increments on this axis reflect orders of magnitude increments in solution time.
The extension with the most significant impact on solution time is the fixed transportation cost (FTC) extension, as indicated by all 24 models without this extension solving the fastest (shown in blue).
Next, we examine how the FTC extension interacts with other extensions. 
The extension that most affects computational time when combined with FTC is the service level agreements (SLA), distinguishing the performance between the green and orange groups. 
Lastly, the addition of  final inventory deviation FID notably impacts solution time as evident in the difference between the orange and red groups.

As expected, any models that include fixed transportation cost (FTC) and service level agreements (SLA either in its simple or window version) require more time (as depicted in red and orange), likely due to the additional binary variables and big-M constraints  \citep{trespalacios2015improved}.
Including the final inventory deviation (FID) extension seems to interact with FTC and SLA, resulting in an increase of computational time. Nevertheless, this difference is reduced as we increase the number of time periods.
In models where FTC is considered and SLAs are not, the computational time decreases independently from all other extensions as shown in green. 
This implies that, on average, FID only represents an increase in solution time if paired with SLA. 
Therefore, if the service level agreement extension is not considered, FID (as well as PATP and NID) can be added without significantly impacting the solution time. 
We conclude that only three different extensions, and their respective interactions, significantly impact the solution time of the model. 
In this sense, the production across time periods (PATP) and negative inventory deviation (NID, which requires additional binary variables) extensions do not have a major impact on solution time.
The model proves to be scalable with time given that even the largest formulations, which yield hundreds of thousands of binary variables, take no more than 9 minutes to solve on average.

\subsection{Disruption Characterization} \label{sec:dis_char}


Effective response optimization and recovery relies on appropriate characterization of the disruption.
For example, consider a scenario where a main reactor within Plant 1 is damaged and needs replacement. 
First, we must evaluate the impact on the plant production  capacity, which involves understanding how production will be affected. 
Relevant information includes: the number of recipes reliant on the damaged reactor, whether there are similar reactors operating in parallel, the existence of alternative routes within the plant can perform the same tasks, and if this disruption implies a plant shut-down.
Second, we need to determine the appropriate time-frame considering factors such as the purchasing lead time for a new reactor, the time required to replace the reactor, and how long it will take to restore the plant to nominal production capacity.
With this information, we can adjust the production bounds to capture the impact of the disruption.

\begin{figure}[!htb]
\centering
\begin{subfigure}{0.49\textwidth}
    \includegraphics[width=1\textwidth]{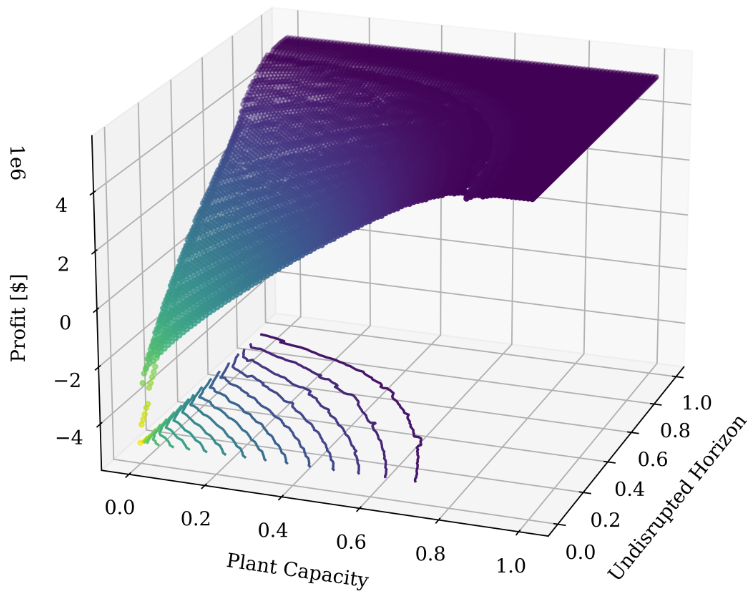}
    \caption{Total profit of the operation.}
    \label{fig:dis_char_obj}
\end{subfigure}
\hfill
\begin{subfigure}{0.49\textwidth}
    \includegraphics[width=1\textwidth]{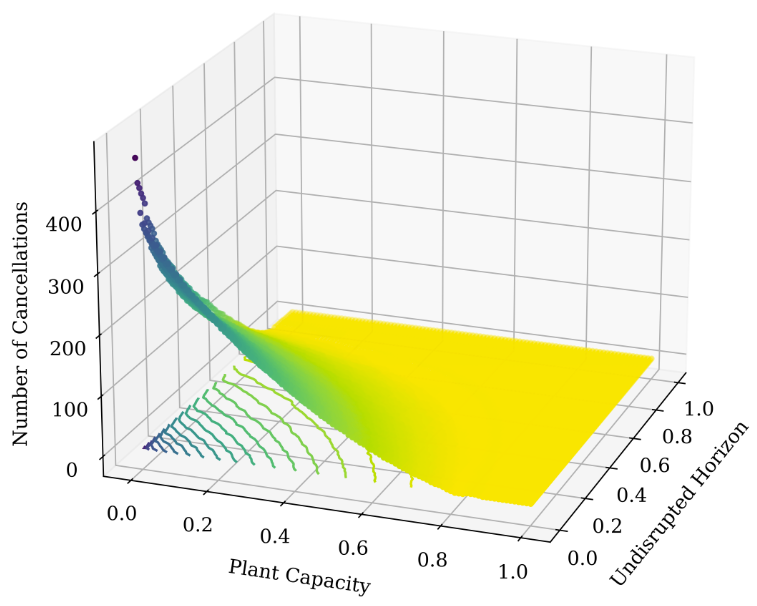}
    \caption{Total number of canceled orders.}
    \label{fig:dis_char_cancel}
\end{subfigure}
        
\caption{Impacts on the profit of the operation and number of orders canceled for different characterizations of a damaged reactor in Plant 1.}
\label{fig:dis_char_objcancel}
\end{figure}

\begin{figure}[!htb]
\centering
\begin{subfigure}{0.49\textwidth}
    \includegraphics[width=1\textwidth]{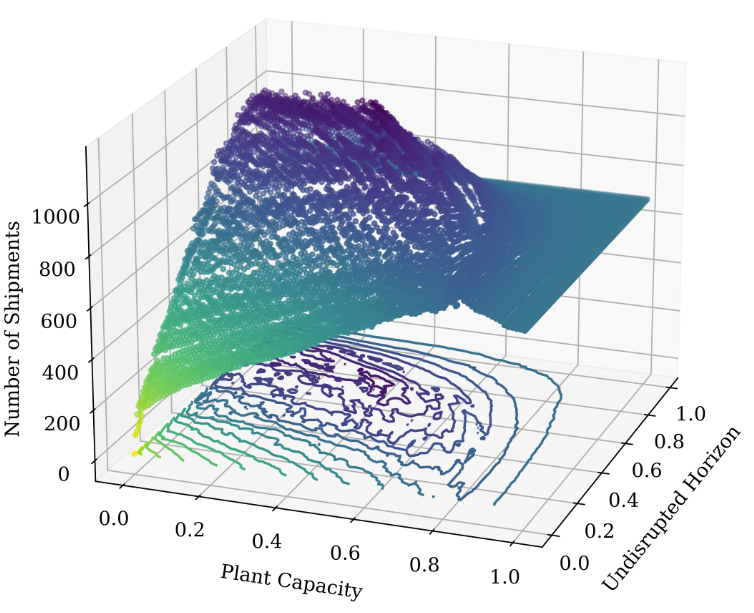}
    \caption{Total number of shipments required in network.}
    \label{fig:dis_char_fin}
\end{subfigure}
\hfill
\begin{subfigure}{0.49\textwidth}
    \includegraphics[width=1\textwidth]{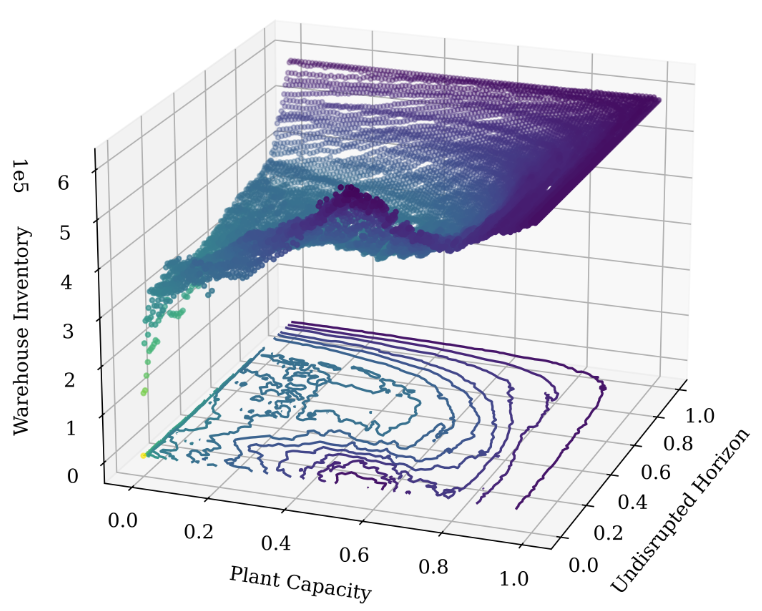}
    \caption{Total amount of inventory level in all warehouses.}
    \label{fig:dis_char_invq}
\end{subfigure}
        
\caption{Impacts on the amount of material flow and warehouse inventory levels for different characterizations of a damaged reactor in Plant 1.}
\label{fig:dis_char_fininv}
\end{figure}

To illustrate the impact of different disruption characterizations, we obtain the optimal disruption response by solving the optimization model for the various disruption scenarios. 
The system was solved over 120 periods with a 12-hour time discretization in a fully reactive context with immediate recovery (as discussed in Section \ref{sec:disruption_modeling}).
Here, we vary plant production bounds and the disruption duration in terms of a fraction.
In both scenarios, a fraction equal to zero yields the most adverse situation, which corresponds to either a complete production shutdown in the plant or a perpetually disrupted state.
Conversely, a fraction of one reflects the nominal state of the system. 
Figure \ref{fig:dis_char_objcancel} offers a view of how these different disruption levels impact both the profit and the number of canceled orders. 
On the other hand, Figure \ref{fig:dis_char_fininv} provides insights into how these disruption levels affect the optimal operation itself.
In this context, we analyze the number of shipments needed and the total inventory in the warehouses of the network. 
It is worth mentioning that these figures represent an overview of the aggregated data for 10 different materials, giving us a big-picture perspective on the solution.
However, the model is still calculating optimal inventory and routing for each material individually, ensuring a tailored product-specific approach to managing supply chain disruptions. 
Note that the color in these figures is proportional to the vertical axis, and it is meant to help visualize the results.


We observe that both subplots depicted in Figure \ref{fig:dis_char_objcancel}, representing the variation in the number of order cancellations and profit, exhibit a monotonic trend, which aligns with our expectation. 
Here, reduced production capacity together with longer disruption periods results in fewer orders being delivered. 
Notably, these behaviors are directly and inversely correlated, as any decision to cancel an order incurs a cost (represented by $\lambda^\delta$) and lowers overall profit.
In addition to illustrating the results of an approach, figures like these are useful for decision-makers to analyze the reliability of their supply chain in the face of disruptions (with optimal decision-making).



The subfigures in Figure \ref{fig:dis_char_fininv} show the number of shipments and total warehouse inventory levels revealing non-monotonic behavior. 
This indicates that different disruption characterizations impact the entire system, not just the plant. The dynamic responses are due to the interconnectivity of the network and complex node relationships, facilitating disruption propagation in time and space.
Additionally, the subfigures exhibit complementarity; within the mid-range of characterization values, there is a region with a high number of shipments and low warehouse inventory levels. This suggests a rapid inventory turnover policy in warehouses as a reactive response to mitigate disruption effects. While the figures provide aggregated high-level information, further investigation of per-scenario models confirmed that high inventory turnover is an optimal response policy within this range.

Both Figures \ref{fig:dis_char_objcancel} and \ref{fig:dis_char_fininv} exhibit a plateau as we approach the undisrupted conditions (e.g., fraction values close to one in both axes).
Here, the model demonstrates the capability to achieve performance levels closely resembling those of an undisrupted operation when either the impact of the disruption on plant production capacity is limited or it does not affect the system for an extended duration.
The behavior makes intuitive sense in Figure \ref{fig:dis_char_objcancel}, where near-undisrupted scenarios consistently yield high profits and eliminate the need for order cancellations and delays.
However, in the operation plots in Figure \ref{fig:dis_char_fininv}, a plateau near the values of near-undisrupted characterization is observed, which is unexpected since these models need not behave similarly in aggregation. Additionally, Figure \ref{fig:dis_char_fin} shows the plateau positioned lower than some middle values, indicating this is not due to variables being close to their bounds but rather reflects genuine similarity in optimal schedules under these scenarios. This demonstrates the  capability of the model to deliver consistent performance under near-undisrupted conditions despite variations in the operational context.


When analyzing Figures \ref{fig:dis_char_objcancel} and \ref{fig:dis_char_fininv}, we observe that except for the origin representing an indefinite plant shutdown, the studied variables change continuously and gradually without abrupt shifts. This consistency holds even in the non-smooth plots in Figure \ref{fig:dis_char_fininv}, indicating that within a narrow range of disruption characterizations the output variable values remain close. Thus, with a good understanding of expected plant capacity and disruption duration, similar aggregated metrics can be expected. While individual schedules may differ slightly within this range, overall patterns exhibit similarity and continuity.

Given the inherent complexity of supply chain and manufacturing networks, accurate characterization of disruptions is a challenging task.
Disruptions can interact, impacting multiple network nodes simultaneously, making it difficult to precisely measure their effects on the entire system. 
However, it is crucial to realize that precise characterization is not a requisite for effectively utilizing the proposed model and obtaining value from it. 
As demonstrated in Section \ref{sec:scalability}, the model is computational efficient, enabling the simultaneous evaluation of a multitude of scenarios.
For instance, in generating the results presented in this section, which encompassed schedules for a two-month operation, more than 10,200 different scenarios were assessed and run in under 30 minutes of computation time. 
The scalability of the model allows decision-makers to utilize it as a tool for optimal ``simulation'', affording the opportunity to experiment with diverse characterizations and evaluate the resulting schedules. 
This feature, coupled with the rolling horizon approach, makes the proposed model as a valuable instrument for conducting ``what-if" scenarios, facilitating responsive decision-making as information is revealed with time and disruptions become better characterized.


\subsection{Late Delivery and Cancellation Cost Sensitivity} \label{sec:sensitivty_costs}

As we discussed in Section \ref{sec:obj_and_bounds}, quantifying the impact of order cancellations ($\lambda^\delta$) and late deliveries ($\lambda^U$) is not necessarily straightforward.
Beyond the costs alone, these can have an impact on the reputation of the company and customer relations.  
In this section we conducted a sensitivity analysis that provides insights into how cancellation and late deliveries affect the model, not only in terms of the primary objective function, but also in terms of the variables that affect the goodwill of the company.


Consider a case where an unplanned disruption occurs due to road maintenance, entirely blocking the route from Plant 1 to Warehouse 1 for a full month. 
In this scenario, the decision-maker is tasked with determining the optimal reactive response for the operations for the next two months in a 12-hour time discretization. 
Achieving this schedule requires the use of alternative routing, like airfreight, and a smart allocation of inventory downstream to overcome the trucking limitation.

\begin{figure}[!htb]
\centering
\begin{subfigure}{0.49\textwidth}
    \includegraphics[width=1\textwidth]{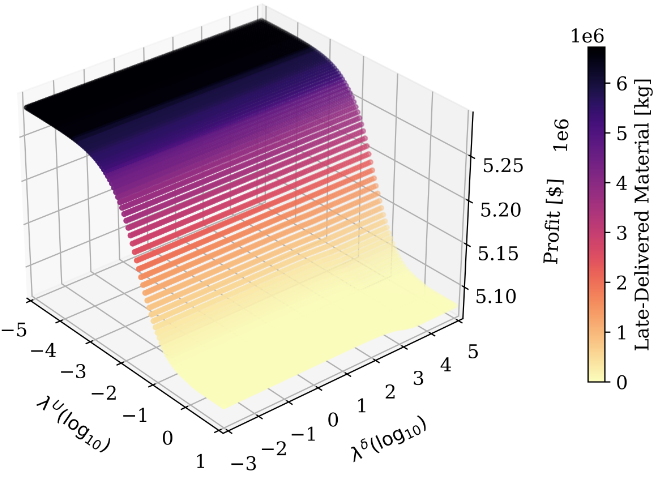}
    \caption{Total profit of the operation.}
    \label{fig:cost_sens_double_late}
\end{subfigure}
\hfill
\begin{subfigure}{0.49\textwidth}
    \includegraphics[width=1\textwidth]{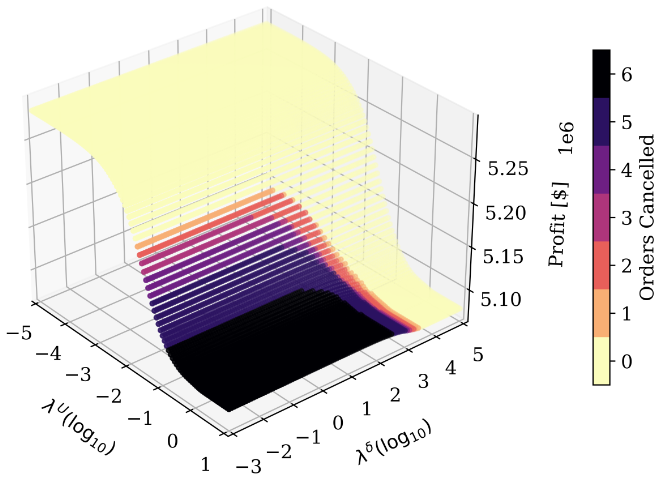}
    \caption{Total number of canceled orders.}
    \label{fig:cost_sens_double_cancels}
\end{subfigure}
        
\caption{Impacts on the profit of the operation and number of orders canceled for different characterizations of a damaged reactor in Plant 1.}
\label{fig:cost_sens_double_latecancels}
\end{figure}

Figures \ref{fig:cost_sens_double_late} and \ref{fig:cost_sens_double_cancels} provide an overview of the effects that varying late delivery and cancellation costs have on the overall profit for this system.
Figure \ref{fig:cost_sens_double_late} shows the total quantity of materials delivered late, and Figure \ref{fig:cost_sens_double_cancels} shows the number of orders that had to be canceled.
We select the same cost for all materials, costumers, and periods to simplify the analysis, and all materials are aggregated to give a holistic overview of the system.

The figures show that the profit of the operation is primarily influenced by the late delivery cost, with the cancellation cost having limited impact.
More specifically, while profit values remain relatively constant across various cancellation penalties, there is a distinct reduction in profit as the cost associated with late deliveries increases.
Notably, the values become stable on both ends of the $\lambda^U$ spectrum implying that the range of costs has been sufficiently explored.
For low values, the model initially refrained from delivering any of the orders affected by the disruption because delivering late was inexpensive. 
On the opposite end, the model decided not to make any late deliveries, even when facing the highest cost penalties, as it had preferred to cancel them.

Figure \ref{fig:cost_sens_double_late} confirms our intuition about covering the full spectrum of late delivery costs based on the color gradient. The top section is uniformly black, indicating that, for a range of late delivery costs, all affected orders were consistently delivered late. 
Here, operational adjustments like re-routing, air freight, and production schedule changes outweigh the income from on-time deliveries. The low late delivery penalty leads the model to opt for paying the accumulated late delivery cost rather than delivering on time.
In contrast, the lower portion of the plot is flat and light yellow, indicating no late deliveries for these cost values. The model prioritizes cost-effective late deliveries over on-time deliveries based on the trade-off between operational adjustments and associated income.

In the upper part of Figure \ref{fig:cost_sens_double_cancels}, the model refrains from canceling orders since delays are not as costly. 
As delays and cancellations become closer in cost, the model identifies a cost optimal decision that is impacted by both parameters. 
Overall, these behaviors align with our expectations considering the nature of the cost penalties. 
With late deliveries being penalized on a per-unit basis, the model retains better control over their allocation. Conversely, the cost of cancellation remains fixed, regardless of the order size, which restricts the flexibility in managing these penalties. This key difference in cost structures influences the strategic choices when dealing with late deliveries and order cancellations in response to supply chain disruptions.

Once the model is solved, it is crucial to evaluate the resulting fill rate within the supply chain, as demonstrated  using the color representations. The reason behind this lies in the fact that customer satisfaction is not always precisely encapsulated within the purely economic objective. 
Instead, we seek a balanced consideration of the trade-off between customer satisfaction and the overall cost of the supply chain.

Conducting a sensitivity analysis can be a useful tool for fine-tuning the cost penalties in order to find the right balance between stressing the system and maintaining customer satisfaction at the desired level. 
Fortunately, these analyses can be performed rapidly thanks to the demonstrated 
solution time, which allows for the quick exploration of multiple scenarios in a manageable and practical manner. 
This adaptive approach ensures that the supply chain is well-prepared to respond to a spectrum of disruptions with a tailored cost strategy.







\section{Conclusions and Future Work} \label{sec:conclusions}

In this study, we present a multiperiod mixed-integer linear programming model designed to optimize the operational response of a supply chain and manufacturing network under disruption. 
The model generates integrated production, shipping, and procurement schedules to mitigate the economic impacts of disruptions. 
It also incorporates customer management decisions such as on-time or delayed deliveries, the extent of delays, and order cancellations.
We extend the model to account for practical supply chain dynamics including time-dependent production, service level agreements, fixed transportation costs, and inventory deviations. Our analysis shows that these extensions have minimal impact on solution time, while identifying which factors and their interactions most affect computational performance.

The model is flexible, allowing for independent or combined disruption scenarios. 
It bridges gaps in the literature by addressing unexpected demands as disruptions and offering effective solutions. 
Using a real-life case study from the chemical industry, we demonstrate the ability of the model to provide tractable solutions as the discretization scheme is refined. 
We also explore the effects of various disruption scenarios and cost penalties for late deliveries and order cancellations showing that, despite significant differences in these factors, the model remains computationally efficient enabling robust ``what-if" analyses to guide decision-making.

As a future research direction, we aim to investigate the scalability of the model concerning network size and assess whether different configurations of arbitrary networks yield varied computational performance.
Currently, the problem is addressed deterministically, with uncertainty indirectly addressed through multiple deterministic instances in a simulation scheme, capitalizing on the computational efficiency of the model.
However, future research will focus on directly incorporating uncertainty by modeling it in the disruption characterization.
Moreover, chance constraints will be explored to ensure the enforcement of a minimum service level within this framework.

\section*{Acknowledgments and Disclaimers} \label{sec:ack_and_dis}

We gratefully acknowledge the Dow Chemical Company for motivating this work and providing a real-life case study. 
Parts of this work were presented in the proceedings of the FOCAPO/CPC 2023 Conference.

\section*{Funding}
We thank the Center for Advanced Process Decision-making (CAPD) and the Chemical Engineering Department of Carnegie Mellon University for funding this work.

\pagebreak
\bibliography{bibliography.bib}

\pagebreak
\appendix
\section{Appendix: Nomenclature} \label{sec:nomenclature}
This section provides a comprehensive nomenclature of all sets, parameters, and variables used in the base model and its proposed extensions.
\subsection{Sets}
\begin{itemize}
    \item Nodes ($\mathcal{N}:=\mathcal{S}\cup\mathcal{P}\cup\mathcal{W}\cup\mathcal{C}$)
    \begin{itemize}
        \item Suppliers ($\mathcal{S}$)
        \item Plants ($\mathcal{P}$)
        \item Warehouses ($\mathcal{W}$)
        \item Customers ($\mathcal{C}$)
        \item Suppliers that require service level agreements ($\mathcal{S}^{SLA} \subseteq \mathcal{S}$)
    \end{itemize}
    \item Arcs ($\mathcal{A}:= (\mathcal{N} \times \mathcal{N}))$ 
    \begin{itemize}
         \item Arcs coming into node $n \in \mathcal{N}$ ($\mathcal{A}^{In}_n := \{(i,j)\in \mathcal A: j=n \}$)
        \item Arcs coming out of node $n \in \mathcal{N}$ ($\mathcal{A}^{Out}_n := \{(i,j)\in \mathcal A: i=n \}$)
    \end{itemize}
    \item Materials ($\mathcal{M}$)
    \begin{itemize}
        \item Materials that can exist at node $n \in \mathcal{N}$ ($\mathcal{M}_n$)
        \item Materials that can move through at arc $a \in \mathcal{A}$ ($\mathcal{M}_a$)
    \end{itemize}
    \item Time periods ($\mathcal{T}$)
    \item Recipes that occur in plant $p \in \mathcal{P}$ ($\mathcal{R}_p$)
\end{itemize}

\subsection{Continuous Variables ($\in \mathbb{R}^+_*$)}
\begin{itemize}
    \item $F^{In}_{mat}$: Flow of material $m \in \mathcal{M}_a$ that is going into arc $a \in \mathcal{A}$ during time period $t \in \mathcal{T}$.
    \item $F^{Out}_{mat}$: Flow of material $m \in \mathcal{M}_a$ that is going out of arc $a \in \mathcal{A}$ during time period $t \in \mathcal{T}$.
    \item $P_{prt}$: Production of recipe $r \in \mathcal{R}_p$ in plant $p \in \mathcal{P}$ during time period $t \in \mathcal{T}$.
    \item $I_{mnt}$: Inventory of material $m \in \mathcal{M}_n$ in node $n \in \mathcal{P} \cup \mathcal{W}$ during time period $t \in \mathcal{T}$.
    \item $B_{mst}$: Amount to buy of material $m \in \mathcal{M}_s$ from supplier $s \in \mathcal{S}$ during time period $t \in \mathcal{T}$.
    \item $D_{mct}$: Demand met of material $m \in \mathcal{M}_c$ to customer $c \in \mathcal{C}$ during time period $t \in \mathcal{T}$.
    \item $U_{mct}$: Unmet demand met of material $m \in \mathcal{M}_c$ to customer $c \in \mathcal{C}$ during time period $t \in \mathcal{T}$.
    \item $P^{In}_{prt}$: Production of recipe $r \in \mathcal{R}_p$ in plant $p \in \mathcal{P}$ starting at time period $t \in \mathcal{T}$.
    \item $P^{Out}_{prt}$: Production of recipe $r \in \mathcal{R}_p$ in plant $p \in \mathcal{P}$ concluding at time period $t \in \mathcal{T}$.
    \item $Deviation_{mn}$: Absolute value of difference between the initial and final inventories of material $m \in \mathcal{M}_n$ in node $n \in \mathcal{P} \cup \mathcal{W}$.
    \item $K_{mnt}$: Distance between the inventory of material $m \in \mathcal{M}_n$ in node $n \in \mathcal{P} \cup \mathcal{W}$ and the threshold value $I^L_{mnt}$ at time period $t \in \mathcal{T}$. This distance is also referred in the paper as negative inventory deviation.
\end{itemize}

\subsection{Binary Variables ($\in \{0,1\}$)}
\begin{itemize}
    \item $y_{mct}$: Indicates if the order of material $m \in \mathcal{M}_c$ to customer $c \in \mathcal{C}$ during time period $t \in \mathcal{T}$ is canceled ($y_{mct}=1$) or not ($y_{mct}=0$).
    \item $x_{mat}$: Indicates the existence of flow of material $m \in \mathcal{M}_a$ going into arc $a \in \mathcal{A}$ during time period $t \in \mathcal{T}$.
    \item $w_{mst}$: Determines if the service level agreement of material $m \in \mathcal{M}_s$ and supplier $s \in \mathcal{S}^{SLA}$ is executed at time period $t \in \mathcal{T}$ or not.
    \item $z_{mnt}$: Indicates that the inventory of material $m \in \mathcal{M}_n$ in node $n \in \mathcal{P} \cup \mathcal{W}$ is below  $I^L_{mnt}$ at time period $t \in \mathcal{T}$.
\end{itemize}

\subsection{Parameters}
\begin{itemize}
    \item $\phi_{rm}$: Mass production (+) or consumption (-) of material $m \in \mathcal{M}_p$ in for recipe $r \in \mathcal{R}_p$ in plant $p \in \mathcal{P}$ during time period $t \in \mathcal{T}$.
    \item $\delta_{mct}$: Order of material $m \in \mathcal{M}_c$ placed by customer $c \in \mathcal{C}$ during time period $t \in \mathcal{T}$.
    \item $SS_{mnt}$: Buffer stock of material $m \in \mathcal{M}_n$ in node $n \in \mathcal{P} \cup \mathcal{W}$ during time period $t \in \mathcal{T}$.
    \item $\alpha_{mnt}$: Fraction of the buffer stock of material $m \in \mathcal{M}_n$ in node $n \in \mathcal{P} \cup \mathcal{W}$ during time period $t \in \mathcal{T}$ that should be enforced regardless of the disruption.
    \item $B^{SLA}_{mst}$: Minimum purchase amount stated by the service level agreement of material $m \in \mathcal{M}_s$ and supplier $s \in \mathcal{S}^{SLA}$ executed at time period $t \in \mathcal{T}$.
    \end{itemize}

    \subsubsection{Variable bounds}
    \begin{itemize}
        \item $\left[F^{L}_{mat}, F^{U}_{mat}\right]$: Lower and upper bounds for the flow of material $m \in \mathcal{M}_a$ that traversing arc $a \in \mathcal{A}$ during time period $t \in \mathcal{T}$.
        \item $\left[P^{L}_{prt}, P^{U}_{prt}\right]$: Lower and upper bounds for the production of recipe $r \in \mathcal{R}_p$ in plant $p \in \mathcal{P}$ during time period $t \in \mathcal{T}$.
        \item $\left[I^{L}_{mnt}, I^{U}_{mnt}\right]$: Lower and upper bounds for the inventory of material $m \in \mathcal{M}_n$ in node $n \in \mathcal{P} \cup \mathcal{W}$ during time period $t \in \mathcal{T}$.
        \item $\left[B^{L}_{mst}, B^{U}_{mst}\right]$: Lower and upper bounds for the amount to buy of material $m \in \mathcal{M}_s$ from supplier $s \in \mathcal{S}$ during time period $t \in \mathcal{T}$.
        \item $\left[D^{L}_{mct}, D^{U}_{mct}\right]$: Lower and upper bounds for the demand met of material $m \in \mathcal{M}_c$ to customer $c \in \mathcal{C}$ during time period $t \in \mathcal{T}$.
        \item $\left[U^{L}_{mct}, U^{U}_{mct}\right]$: Lower and upper bounds for the unmet demand met of material $m \in \mathcal{M}_c$ to customer $c \in \mathcal{C}$ during time period $t \in \mathcal{T}$.
    \end{itemize}
    
    \subsubsection{Initial conditions}
    \begin{itemize}
        \item $F^{In,0}_{ma}$: Pre-planned flow of material $m \in \mathcal{M}_a$ that is going into arc $a \in \mathcal{A}$ during time period $t \in \mathcal{T}$ in absence of disruption.
        \item $F^{Out,0}_{ma}$: Pre-planned flow of material $m \in \mathcal{M}_a$ that is going out of arc $a \in \mathcal{A}$ during time period $t \in \mathcal{T}$ in absence of disruption.
        \item $P^0_{pr}$:  Pre-planned production of recipe $r \in \mathcal{R}_p$ in plant $p \in \mathcal{P}$ during time period $t \in \mathcal{T}$ in absence of disruption.
        \item $I^0_{mn}$: Existent inventory of material $m \in \mathcal{M}_n$ in node $n \in \mathcal{P} \cup \mathcal{W}$ during time period $t \in \mathcal{T}$ at the time disruption occurs.
        \item $B^0_{ms}$:  Pre-planned amount to buy of material $m \in \mathcal{M}_s$ from supplier $s \in \mathcal{S}$ during time period $t \in \mathcal{T}$ in absence of disruption.
        \item $D^0_{mc}$:  Pre-planned demand met of material $m \in \mathcal{M}_c$ to customer $c \in \mathcal{C}$ during time period $t \in \mathcal{T}$ in absence of disruption.
    \end{itemize}
    
    \subsubsection{Time durations ($\tau$)}
    \begin{itemize}
        \item $\tau_{ma{t}}$: Time it takes material $m \in \mathcal{M}_a$ to traverse arc $a \in \mathcal{A}$ when shipped at time period $t \in \mathcal{T}$.
        \item $\tau^P_{pr{t}}$: Time it takes to execute recipe $r \in \mathcal{R}_p$ in plant $p \in \mathcal{P}$ starting at time period $t \in \mathcal{T}$.
        \item $\tau^{SLA}_{ms{t}}$: Time of the window to fulfill a service level agreement of material $m \in \mathcal{M}_s$ to supplier $s \in \mathcal{S}^{SLA}$ for a contract that started at time period $t \in \mathcal{T}$.
    \end{itemize}
    
    \subsubsection{Economic parameters ($\lambda$)}
    \begin{itemize}
        \item $\lambda^D_{mct}$: Price per unit of material $m \in \mathcal{M}_c$ sold to customer $c \in \mathcal{C}$ during time period $t \in \mathcal{T}$.
        \item $\lambda^U_{mct}$: Cost penalty per unit of unmet demand of material $m \in \mathcal{M}_c$ to customer $c \in \mathcal{C}$ during time period $t \in \mathcal{T}$.
        \item $\lambda^\delta_{mct}$: Cost penalty of canceling an order of material $m \in \mathcal{M}_c$ to customer $c \in \mathcal{C}$ during time period $t \in \mathcal{T}$.
        \item $\lambda^F_{mat}$: Cost of transporting a unit material $m \in \mathcal{M}_a$ through arc $a \in \mathcal{A}$ during time period $t \in \mathcal{T}$.
        \item $\lambda^B_{mst}$: Cost of buying a unit material $m \in \mathcal{M}_s$ from supplier $s \in \mathcal{S}$ during time period $t \in \mathcal{T}$.
        \item $\lambda^P_{prt}$: Cost of producing a unit of recipe $r \in \mathcal{R}_p$ in plant $p \in \mathcal{P}$ during time period $t \in \mathcal{T}$.
        \item $\lambda^I_{mnt}$: Cost of holding a unit material $m \in \mathcal{M}_n$ as inventory in node $n \in \mathcal{P} \cup \mathcal{W}$ during time period $t \in \mathcal{T}$.
         \item $\lambda^{Ffix}_{mat}$: Fixed cost of transporting material $m \in \mathcal{M}_a$ through arc $a \in \mathcal{A}$ during time period $t \in \mathcal{T}$.
         \item $\lambda^{Dev}_{mn}$: Penalty for the deviation between the initial and final inventories of material $m \in \mathcal{M}_n$ in node $n \in \mathcal{P} \cup \mathcal{W}$.
         \item $\lambda^K_{mnt}$: Penalty for the deviation between the inventory of material $m \in \mathcal{M}_n$  in node $n \in \mathcal{P} \cup \mathcal{W}$ and the threshold value $I^L_{mnt}$ during time period $t \in \mathcal{T}$.

    \end{itemize}

\end{document}